\documentclass[10pt]{article}

\RequirePackage{amsmath}
\RequirePackage{amsfonts}
\RequirePackage{amssymb}
\RequirePackage{hyperref}
\RequirePackage{multirow}
\RequirePackage{algorithm}
\RequirePackage{algpseudocode}
\RequirePackage{graphicx}
\RequirePackage{xcolor}
\RequirePackage{amsthm}
\RequirePackage{amsmath}
\RequirePackage{amsfonts}
\RequirePackage{amssymb}
\RequirePackage{hyperref}
\RequirePackage{multirow}
\RequirePackage{algorithm}
\RequirePackage{algpseudocode}
\RequirePackage{graphicx}
\RequirePackage{xcolor}
\allowdisplaybreaks

\usepackage{dpr_fec,ifthen,dsfont}

\newtheorem{theorem}{Theorem}

\newtheorem{example}[theorem]{Example}

\newtheorem{proposition}[theorem]{Proposition}
\newtheorem{remark}[theorem]{Remark}

\usepackage{isetechreport}
\def\reportyear{20T}
\def\reportno{012}
\def\revisionno{0}
\def\originaldate{\today}
\coraltrue
\cvcrfalse

\begin{document}

\title {Sequential Quadratic Optimization for Nonlinear Equality Constrained Stochastic Optimization}

\author{Albert Berahas\thanks{E-mail: albertberahas@gmail.com}}
\affil{Department of Industrial and Operations Engineering, University of Michigan}
\author{Frank E.~Curtis\thanks{E-mail: frank.e.curtis@gmail.com}}
\author{Daniel P.~Robinson\thanks{E-mail: daniel.p.robinson@gmail.com}}
\author{Baoyu Zhou\thanks{E-mail: baz216@lehigh.edu}}
\affil{Department of Industrial and Systems Engineering, Lehigh University}

\titlepage

\maketitle

\begin{abstract}
  Sequential quadratic optimization algorithms are proposed for solving smooth nonlinear optimization problems with equality constraints.  The main focus is an algorithm proposed for the case when the constraint functions are deterministic, and constraint function and derivative values can be computed explicitly, but the objective function is stochastic.  It is assumed in this setting that it is intractable to compute objective function and derivative values explicitly, although one can compute stochastic function and gradient estimates.  As a starting point for this stochastic setting, an algorithm is proposed for the deterministic setting that is modeled after a state-of-the-art line-search SQP algorithm, but uses a stepsize selection scheme based on Lipschitz constants (or adaptively estimated Lipschitz constants) in place of the line search.  This sets the stage for the proposed algorithm for the stochastic setting, for which it is assumed that line searches would be intractable.  Under reasonable assumptions, convergence (resp.,~convergence in expectation) from remote starting points is proved for the proposed deterministic (resp.,~stochastic) algorithm.  The results of numerical experiments demonstrate the practical performance of our proposed techniques.
\end{abstract}


\section{Introduction}\label{sec.introduction}

We consider the design of algorithms for solving smooth nonlinear optimization problems with equality constraints.  Such problems arise in various important applications throughout science and engineering.

Numerous algorithms have been proposed for solving \emph{deterministic} equality constrained optimization problems.  Penalty methods \cite{Cour43,Flet87}, including augmented Lagrangian methods \cite{ConnGoulToin92,Hest69,Powe69}, attempt to solve such problems by penalizing constraint violation through an objective term---weighted by a penalty parameter---and employing unconstrained optimization techniques for solving (approximately) a corresponding sequence of penalty subproblems.  Such algorithms can behave poorly due to ill-conditioning and/or nonsmoothness of the penalty subproblems, depending on the type of penalty function employed.  Their performance also often suffers due to sensitivity to the scheme for updating the penalty parameter.

Algorithms that consistently outperform penalty methods are those based on sequential quadratic optimization (commonly known as SQP), which in this setting of equality constrained optimization is intimately connected to the idea of applying Newton's method to stationarity conditions of the problem \cite{Wils63}.  In particular, it is commonly accepted that one of the state-of-the-art algorithms for solving equality constrained optimization problems is such an SQP method that chooses stepsizes based on a line search applied to an exact penalty function \cite{Han77,HanMang79,Powe78}.  In such a method, the penalty function acts as a merit function only. It does not influence the computed search direction; it only influences the computed stepsize.

Significantly fewer algorithms have been proposed for solving \emph{stochastic} equality constrained optimization problems.  In particular, in this paper, we focus on such problems with constraint functions that are deterministic, but objective functions that are stochastic, in the sense that the objective is an expectation of a function defined with respect to a random variable with unknown distribution.  (Various modeling paradigms have been proposed for solving problems involving stochastic constraints.  These are out of our scope; we refer the reader to~\cite{ShapDentRusz09}.)  We assume that it is intractable to compute objective function and gradient values, although one is able to compute (unbiased) stochastic gradient estimates.  A few algorithms have been proposed that may be employed in this setting \cite{ChenTungVeduMori18,KumaSoumMhamHara18,NandPathAbhiSing19,RaviDinhLokhSing19}, although these are based on penalty methodologies, so do not benefit from advantages of SQP techniques.  Let us also mention various proposed stochastic Frank-Wolfe algorithms \cite{goldfarb2017linear,hazan2016variance,locatello2019stochastic,lu2020generalized,negiar2020stochastic,reddi2016stochastic,zhang2020one} for (non)convex stochastic optimization with convex constraints, although these are not applicable for our setting of having general nonlinear equality constraints.

\subsection{Contributions}

In this paper, we propose two algorithms modeled after the aforementioned line-search SQP methodology.  Our primary focus is an algorithm for the aforementioned setting of a problem with deterministic constraint functions, but a stochastic objective function.  However, as a first step for considering this setting, we begin by proposing an algorithm for the deterministic setting that employs an adaptive stepsize selection scheme that makes use of Lipschitz constants (or adaptively updated Lipschitz constant estimates) rather than a line search.  Based on this algorithm for the deterministic setting, we propose our algorithm for the stochastic setting that also uses Lipschitz constants (or, in practice, estimates of them) for stepsize selection.

We prove under common assumptions that our deterministic algorithm has convergence guarantees that match those of a state-of-the-art line-search SQP method.  In addition, we prove under loose assumptions that our stochastic algorithm offers convergence guarantees that can match those of our deterministic algorithm \emph{in expectation}.  In particular, the results that we prove for our stochastic algorithm are of the type offered by stochastic gradient schemes for unconstrained optimization \cite{BottCurtNoce18}.  An additional challenge for constrained stochastic optimization is potentially poor behavior of an adaptive merit function parameter that balances emphasis between minimizing constraint violation and reducing the objective function.  To address this, in addition to our aforementioned convergence analysis, which considers the behavior of the algorithm under good behavior of this adaptive parameter, we prove under pragmatic assumptions that certain poor behavior either cannot occur or only occurs in extreme circumstances, and other poor behavior occurs with probability zero.

The results of numerical experiments show that our deterministic algorithm is as reliable as a state-of-the-art line-search SQP method, although, as should be expected, it is sometimes less efficient than such a method that performs line searches.  Our experiments with our stochastic algorithm show that it consistently and significantly outperforms an approach that attempts to solve constrained problems by applying a stochastic (sub)gradient scheme to minimize an exact penalty function.

\subsection{Notation}

Let $\R{}$ denote the set of real numbers (i.e., scalars), let $\R{}_{\geq r}$ (resp.,~$\R{}_{>r}$) denote the set of real numbers greater than or equal to (resp.,~greater than) $r \in \R{}$, let $\R{n}$ denote the set of $n$-dimensional real vectors, let $\R{m \times n}$ denote the set of $m$-by-$n$-dimensional real matrices, and let $\mathbb{S}^n$ denote the set of $n$-by-$n$-dimensional symmetric matrices.  The set of natural numbers is denoted as $\N{} := \{0,1,2,\dots\}$.  For any $m \in \N{}$, let $[m]$ denote the set of integers $\{1,\dots,m\}$.

Each of our algorithms is iterative, generating a sequence of iterates $\{x_k\}$ with $x_k \in \R{n}$ for all $k \in \N{}$.  The iteration number is also appended as a subscript to other quantities corresponding to each iteration; e.g., $f_k := f(x_k)$ for all $k \in \N{}$.

\subsection{Organization}

Our algorithm for the deterministic setting is proposed and analyzed in \S\ref{sec.deterministic}.  We present our analysis alongside that of a line-search SQP method for ease of comparison with this state-of-the-art strategy.  Our algorithm for the stochastic setting is proposed and analyzed in \S\ref{sec.stochastic}.  The results of numerical experiments are provided in \S\ref{sec.numerical} and concluding remarks are offered in \S\ref{sec.conclusion}.

\section{Deterministic Setting}\label{sec.deterministic}

Given an objective function $f : \R{n} \to \R{}$ and a constraint function $c : \R{n} \to \R{m}$, consider the optimization problem
\bequation\label{prob.deterministic}
  \min_{x\in\R{n}}\ f(x)\ \st\ c(x) = 0.
\eequation
We make the following assumption about the optimization problem~\eqref{prob.deterministic} and the algorithms that we propose, each of which generates an iterate sequence $\{x_k\} \subset \R{n}$, search direction sequence $\{d_k\} \subset \R{n}$, and trial stepsize sequence $\{\alpha_{k,j}\} \subset \R{}_{>0}$.

\bassumption\label{ass.deterministic}
  Let $\Xcal \subseteq \R{n}$ be an open convex set containing the iterates $\{x_k\}$ and trial points $\{x_k+\alpha_{k,j}d_k\}$.  The objective function $f : \R{n} \to \R{}$ is continuously differentiable and bounded below over $\Xcal$, and its gradient $\nabla f : \R{n} \to \R{n}$ is Lipschitz continuous with constant $L$ and bounded over $\Xcal$.  The constraint function $c : \R{n} \to \R{m}$ and its Jacobian $\nabla c^T : \R{n} \to \R{m \times n}$ are bounded over $\Xcal$, each gradient $\nabla c_i : \R{n} \to \R{n}$ is Lipschitz continuous with constant $\gamma_i$ over $\Xcal$ for all $i \in \{1,\dots,m\}$, and the singular values of $\nabla c^T$ are bounded away from zero over $\Xcal$.
\eassumption

Most of the statements in Assumption~\ref{ass.deterministic} are standard smoothness assumptions for the objective and constraint functions.  We remark that we do not assume that the set $\Xcal$ is bounded.  The assumption that the singular values of $\nabla c^T$ are bounded away from zero is equivalent to the linear independence constraint qualification (LICQ).  The LICQ is a relatively strong assumption in the modern literature on algorithms for solving constrained optimization problems, but it is a reasonable one in our context due to the significant challenges that arise in the stochastic setting in \S\ref{sec.stochastic}.

Defining the Lagrangian $\ell : \R{n} \times \R{m} \to \R{}$ corresponding to \eqref{prob.deterministic} by $\ell(x,y) = f(x) + c(x)^Ty$, first-order stationarity conditions for \eqref{prob.deterministic}---which are necessary due to the inclusion of the LICQ in Assumption~\ref{ass.deterministic}---are given by
\bequation\label{eq.KKT}
  0 = \bbmatrix \nabla_x \ell(x,y) \\ \nabla_y \ell(x,y) \ebmatrix = \bbmatrix \nabla f(x) + \nabla c(x)y \\ c(x) \ebmatrix.
\eequation

A consequence of Lipschitz continuity of the constraint functions is the following.  Since this fact is well known and easily proved, we present it without proof.

\blemma\label{lem.c_upper}
  Under Assumption~\ref{ass.deterministic}, it follows for any $x \in \R{n}$, $\alpha \in \R{}_{>0}$, and $d \in \R{n}$ such that $(x,x+\alpha d) \in \Xcal \times \Xcal$ that
  \bequationNN
    \baligned
      |c_i(x + \alpha d)| &\leq |c_i(x) + \alpha \nabla c_i(x)^T d| + \thalf \gamma_i \alpha^2 \|d\|^2_2\ \ \text{for all}\ \ i \in [m] \\ \text{and}\ \ 
      \|c(x + \alpha d)\|_1 &\leq \|c(x) + \alpha \nabla c(x)^T d\|_1 + \thalf \Gamma \alpha^2 \|d\|^2_2\ \ \text{with}\ \ \Gamma := \sum_{i\in[m]} \gamma_i.
    \ealigned
  \eequationNN
\elemma

\subsection{Merit Function}\label{sec.Merit}

As is common in SQP techniques, our algorithms use as a merit function the $\ell_1$-norm penalty function $\phi : \R{n} \times \R{}_{>0} \to \R{}$ defined by
\bequation\label{eq.penalty_function}
  \phi(x,\tau) = \tau f(x) + \|c(x)\|_1.
\eequation
Here, $\tau \in \R{}_{>0}$ is a merit parameter, the value of which is chosen in the algorithm according to a positive sequence $\{\tau_k\}$ that is set adaptively.  We make use of a local model of the merit function $q : \R{n} \times \R{}_{>0} \times \R{n} \times \mathbb{S}^n \times \R{n} \to \R{}$ defined by
\bequationNN
  q(x,\tau,\nabla f(x),H,d) = \tau (f(x) + \nabla f(x)^Td + \thalf \max\{d^THd,0\}) + \|c(x) + \nabla c(x)^Td\|_1.
\eequationNN
A critical quantity in our algorithms is the reduction in this model for a given $d \in \R{n}$ with $c(x) + \nabla c(x)^Td = 0$, i.e., $\Delta q : \R{n} \times \R{}_{>0} \times \R{n} \times \mathbb{S}^n \times \R{n} \to \R{}$ defined by
\bequation\label{def.merit_model_reduction}
  \baligned
    \Delta q(x,\tau,\nabla f(x),H,d)
      :=&\ q(x,\tau,\nabla f(x),H,0) - q(x,\tau,\nabla f(x),H,d) \\
       =&\ -\tau(\nabla f(x)^Td + \thalf \max\{d^THd, 0\}) + \|c(x)\|_1.
  \ealigned
\eequation
The following lemma shows an important relationship between the directional derivative of the merit function and this model reduction function.

\blemma\label{lem.directional_derivative}
  Given $(x,\tau,H,d) \in \R{n} \times \R{}_{>0} \times \mathbb{S}^n \times \R{n}$ with $c(x) + \nabla c(x)^Td = 0$,
  \bequation\label{eq.directional_derivative}
    \phi'(x,\tau,d) = \tau \nabla f(x)^Td - \|c(x)\|_1 \leq -\Delta q(x,\tau,\nabla f(x),H,d).
  \eequation
  where $\phi' : \R{n} \times \R{}_{>0} \times \R{n} \to \R{}$ is the directional derivative of $\phi$ at $(x,\tau)$ for $d$.
\elemma
\bproof
  The first equation in \eqref{eq.directional_derivative} is well known; see, e.g., \cite[Theorem~18.2]{NoceWrig06}.  On the other hand, from the definition \eqref{def.merit_model_reduction} one finds that
  \bequationNN
    \Delta q(x,\tau,\nabla f(x),H,d) = -\phi'(x,\tau,d) - \thalf \tau \max\{d^THd,0\} \leq -\phi'(x,\tau,d),
  \eequationNN
  which shows the inequality in \eqref{eq.directional_derivative}.
\eproof

\subsection{Algorithm Preliminaries}

The algorithms that we discuss for solving \eqref{prob.deterministic} are based on an SQP paradigm.  Specifically, at $x_k$ for all $k \in \N{}$, a search direction $d_k \in \R{n}$ is computed by solving a quadratic optimization subproblem based on a local quadratic model of $f$ and a local affine model of $c$ about $x_k$.  Letting $f_k := f(x_k)$, $g_k := \nabla f(x_k)$, $c_k := c(x_k)$, and $J_k := \nabla c(x_k)^T$ for all $k \in \N{}$ and given a sequence~$\{H_k\}$ satisfying Assumption~\ref{ass.H} below (a standard type of sufficiency condition for equality constrained optimization), this subproblem is given by
\bequationNN
  \min_{d\in\R{n}}\ f_k + g_k^Td + \thalf d^TH_kd\ \ \st\ \ c_k + J_kd = 0.
\eequationNN
The optimal solution $d_k$ of this subproblem, and an associated Langrange multiplier $y_k \in \R{m}$, can be obtained by solving the linear system of equations
\bequation\label{eq.system_deterministic}
  \bbmatrix H_k & J_k^T \\ J_k & 0 \ebmatrix \bbmatrix d_k \\ y_k \ebmatrix = - \bbmatrix g_k \\ c_k \ebmatrix.
\eequation

\bassumption\label{ass.H}
  The sequence $\{H_k\}$ is bounded in norm by $\kappa_H \in \R{}_{>0}$.  In addition, there exists a constant $\zeta \in \R{}_{>0}$ such that, for all $k \in \N{}$, the matrix $H_k$ has the property that $u^TH_ku \geq \zeta \|u\|_2^2$ for all $u \in \R{n}$ such that $J_k u = 0$.
\eassumption

We stress that our algorithms and analysis do \emph{not} assume that $H_k$ is equal to the Hessian of the Lagrangian at $x_k$ for some multiplier $y_k$, although choosing $\{H_k\}$ in this manner would be appropriate in order to ensure fast local convergence guarantees.  Since our focus is only on achieving convergence to stationarity from remote starting points, we merely assume that $\{H_k\}$ satisfies Assumption~\ref{ass.H}.

Under Assumptions~\ref{ass.deterministic} and \ref{ass.H}, the following results are well known in the literature.

\blemma\label{lem.nonsingular}
  For all $k \in \N{}$, the linear system~\eqref{eq.system_deterministic} has a unique solution.
\elemma

\blemma\label{lem.stationary}
  For any $k \in \N{}$, the solution $(d_k,y_k)$ obtained by solving~\eqref{eq.system_deterministic} has $d_k = 0$ if and only if the pair $(x_k,y_k)$ satisfies \eqref{eq.KKT}.
\elemma

\subsection{Algorithms}\label{sec.Algorithms}

In this section, we present two algorithms for solving problem~\eqref{prob.deterministic}.  The first algorithm chooses stepsizes based on a rule using Lipschitz constant estimates, which can be set adaptively.  This algorithm is new to the literature and establishes a foundation upon which our method for the stochastic setting will be built.  The second algorithm, by contrast, employs a standard type of backtracking line search.  This algorithm is standard in the literature.  We prove a convergence theory for it alongside that for our newly proposed algorithm for illustrative purposes.

In both algorithms, after $d_k$ is computed, the merit parameter $\tau_k$ is set.  This is done by first setting, for some $\sigma \in (0,1)$, a trial value $\tau_k^{trial} \in \R{}_{>0} \cup \{\infty\}$ by
\bequation\label{eq.merit_parameter_trial}
  \tau_k^{trial} \gets \bcases \infty & \text{if $g_k^Td_k + \max\{d_k^TH_kd_k,0\} \leq 0$} \\ \tfrac{(1 - \sigma)\|c_k\|_1}{g_k^Td_k + \max\{d_k^TH_kd_k,0\}} & \text{otherwise.} \ecases
\eequation
(If $c_k = 0$, then it follows from \eqref{eq.system_deterministic} and Assumption~\ref{ass.H} that $d_k^TH_kd_k \geq 0$ and $g_k^Td_k + d_k^TH_kd_k = 0$, meaning $\tau_k^{trial} \gets \infty$.  Hence, $\tau_k^{trial} < \infty$ requires $\|c_k\|_1 > 0$, in which case $\tau_k^{trial} > 0$.)  Then, the merit parameter $\tau_k$ is set, for some $\epsilon \in (0,1)$, by
\bequation\label{eq.merit_parameter_lower}
  \tau_k \gets \bcases \tau_{k-1} & \text{if $\tau_{k-1} \leq \tau_k^{trial}$} \\ (1-\epsilon) \tau_k^{trial} & \text{otherwise.} \ecases
\eequation
This ensures that $\tau_k \leq \tau_k^{trial}$.  Regardless of the case in \eqref{eq.merit_parameter_lower}, it follows that
\bequation\label{eq.merit_model_reduction_lower}
  \Delta q(x_k,\tau_k,g_k,H_k,d_k) \geq \thalf \tau_k \max\{d_k^TH_kd_k,0\} + \sigma \|c_k\|_1.
\eequation
This inequality will be central in our analysis of both algorithms.  In particular, it will be useful when combined with the fact that each algorithm ensures that, for all $k \in \N{}$, the stepsize $\alpha_k \in \R{}_{>0}$ is selected such that for $\eta \in (0,1)$ one finds
\bequation\label{eq.sufficient_decrease}
  \phi(x_k+\alpha_kd_k,\tau_k) \leq \phi(x_k,\tau_k) - \eta \alpha_k \Delta q(x_k,\tau_k,g_k,H_k,d_k).
\eequation

\begin{remark}
  An alternative approach for setting the merit function parameter, which is commonly found in textbooks on nonlinear constrained optimization, is to set it based on the computed Lagrange multiplier estimate $y_k$.  For example, in the context of our $\ell_1$-norm exact penalty function $\phi(x_k,\cdot)$, one can ensure that the computed search direction $d_k$ is a direction of descent for $\phi(\cdot,\tau_k)$ from $x_k$ if $\tau_k < \|y_k\|_{\infty}^{-1}$; see, e.g., \cite{NoceWrig06}. However, it has been recognized that it is often better in practice to set it based on ensuring sufficient reduction in a model of the merit function (see, e.g., \cite{ByrdGilbNoce00,ByrdHribNoce99}), which is our motivation for using the rule defined by  \eqref{eq.merit_parameter_trial}--\eqref{eq.merit_parameter_lower}.
\end{remark}

Our first algorithm is stated as Algorithm~\ref{alg.sqp_adaptive}.  A signifying feature of it is the manner in which it can adapt Lipschitz constant estimates, which are used in the stepsize selection scheme.  For any $(k,j) \in \N{} \times \N{}$, if the estimates $L_{k,j}$ and $\{\gamma_{k,i,j}\}_{i=1}^m$ satisfy $L_{k,j} \geq L$ and $\gamma_{k,i,j} \geq \gamma_{i}$ for all $i \in [m]$, then it follows (see~\cite{Nest04} and Lemma~\ref{lem.c_upper}) that for $\alpha_{k,j} \in \R{}_{>0}$ yielding $x_k + \alpha_{k,j}d_k \in \Xcal$ (recall Assumption~\ref{ass.deterministic}) one has
\bsubequations\label{eq.Lipschitz_bounds}
  \begin{align}
    f(x_k + \alpha_{k,j} d_k) &\leq f_k + \alpha_{k,j} g_k^Td_k + \thalf L_{k,j} \alpha_{k,j}^2 \|d_k\|_2^2 \label{eq.Lipschitz_bound_f} \\ \text{and}\ \  
    |c_i(x_k + \alpha_{k,j} d_k)| &\leq |c_i(x_k) + \alpha_{k,j} \nabla c_i(x_k)^T d_k| + \thalf \gamma_{k,i,j} \alpha_{k,j}^2 \|d_k\|_2^2 && \label{eq.Lipschitz_bound_c}
\end{align}
\end{subequations}
for all $i \in [m]$.  If one knows Lipschitz constants for $\nabla f$ and $\{\nabla c_i\}_{i=1}^m$, then one could simply set $L_{k,0}$ and $\gamma_{k,i,0}$ for all $i \in [m]$ to these values for all $k \in \N{}$, in which case the inner \textbf{for} loop would terminate in iteration $j=0$ for all $k \in \N{}$.  However, if such Lipschitz constants are unknown, as is often the case, then the adaptive procedure in Algorithm~\ref{alg.sqp_adaptive} ensures that convergence can be guaranteed, as shown in the next subsection.  For now, we simply prove the following lemma showing that the inner loop of the algorithm is well-posed.  (One could choose a different increase factor $\rho \in \R{}_{>1}$ for each Lipschitz constant estimate; we use a common value of $\rho$ for simplicity.)

\blemma\label{lem.well_defined}
  Under Assumption~\ref{ass.deterministic}, the inner \textbf{for} loop in Algorithm \ref{alg.sqp_adaptive} is well-posed in that for any $k \in \N{}$, it terminates finitely.  In addition, for all $k \in \N{}$,
  \bequation\label{eq.L_rho_bounds}
    \baligned
      L_k \leq L_{\max} &:= \max\left\{L_{-1},\rho L\right\} \\
      \text{and}\ \ \gamma_{k,i} \leq \gamma_{\max,i} &:= \max\left\{\gamma_{-1,i},\rho \gamma_{i} \right\}\ \ \text{for all}\ \ i \in [m].
    \ealigned
  \eequation
\elemma
\bproof
  To derive a contradiction, suppose that for some $k \in \N{}$ the inner \textbf{for} loop does not terminate.  This means that for each iteration of the \textbf{for} loop at least one inequality in \eqref{eq.Lipschitz_bounds} does not hold.  In such a case, the \textbf{for} loop sets $L_{k,j+1}$ (resp.,~$\gamma_{k,i,j+1}$ for some $i \in [m]$) as $\rho > 1$ times $L_{k,j}$ (resp.,~$\gamma_{k,i,j}$ for some $i \in [m]$).  This leads to a contradiction to the fact that if $L_{k,j} \geq L$ and $\gamma_{k,i,j} \geq \gamma_{i}$ for all $i \in [m]$, then \eqref{eq.Lipschitz_bounds} holds.  Finally, \eqref{eq.L_rho_bounds} follows from the initialization of the Lipschitz constant estimates; the fact that if any of these values is ever increased in the \textbf{for} loop, then this occurs by the value being multiplied by $\rho > 1$; and the fact that for all $k \in \N{}$ the algorithm initializes $L_{k,0} \in (0,L_{k-1}]$ and $\gamma_{k,i,0} \in (0,\gamma_{k-1,i}]$ for all $i \in [m]$.
\eproof

\balgorithm[ht]
  \caption{SQP Algorithm with Adaptive Lipschitz Constant Estimates}
  \label{alg.sqp_adaptive}
  \balgorithmic[1]
    \Require $x_0 \in \R{n}$; $\tau_{-1} \in \R{}_{>0}$; $\epsilon \in (0,1)$; $\sigma \in (0,1)$; $\eta \in (0,1)$; $\rho \in \R{}_{>1}$; $L_{-1} \in \R{}_{>0}$; $\gamma_{-1,i} \in \R{}_{>0}$ for all $i \in [m]$
    \For{\textbf{all} $k \in \N{}$}
      \State Compute $(d_k,y_k)$ as the solution of \eqref{eq.system_deterministic}
      \State \textbf{if} $(x_k,y_k)$ satisfies \eqref{eq.KKT} \textbf{then return} $(x_k,y_k)$ \label{step.termination_adaptive}
      \State Set $\tau_k^{trial}$ by \eqref{eq.merit_parameter_trial} and $\tau_k$ by \eqref{eq.merit_parameter_lower}
      \State Initialize $L_{k,0} \in (0,L_{k-1}]$ and $\gamma_{k,i,0} \gets (0,\gamma_{k-1,i}]$ for all $i \in [m]$ \label{step.initialize}
      \For{\textbf{all} $j \in \N{}$} \label{step.loop}
        \State Set
        \bequationNN
          \baligned
            \widehat\alpha_{k,j} &\gets \tfrac{2(1-\eta) \Delta q(x_k,\tau_k,g_k,H_k,d_k)}{(\tau_k L_{k,j} + \sum_{i\in[m]} \gamma_{k,i,j})\|d_k\|_2^2}\ \ \text{and} \\ \widetilde\alpha_{k,j} &\gets \widehat\alpha_{k,j} - \tfrac{4 \|c_k\|_1}{(\tau_k L_{k,j} + \sum_{i\in[m]} \gamma_{k,i,j})\|d_k\|_2^2}
          \ealigned
        \eequationNN
        \State Set
		\bequationNN
		  \alpha_{k,j} \gets
		    \bcases
		      \widehat\alpha_{k,j} & \text{if $\widehat\alpha_{k,j} < 1$} \\
		      1 & \text{if $\widetilde\alpha_{k,j} \leq 1 \leq \widehat\alpha_{k,j}$} \\
		      \widetilde\alpha_{k,j} & \text{if $\widetilde\alpha_{k,j} > 1$}
		    \ecases
		\eequationNN
		\State \textbf{if} \eqref{eq.sufficient_decrease} or \eqref{eq.Lipschitz_bounds} holds \textbf{then}
		\State \hspace{\algorithmicindent} Set $L_k \gets L_{k,j}$ and $\gamma_{k,i} \gets \gamma_{k,i,j}$ for all $i \in [m]$
		\State \hspace{\algorithmicindent} Set $\alpha_k \gets \alpha_{k,j}$ and $x_{k+1} \gets x_k + \alpha_kd_k$ and \textbf{break} (loop over $j \in \N{}$) \label{step.term_loop}
		\State \textbf{else}
		\State \hspace{\algorithmicindent} \textbf{if} \eqref{eq.Lipschitz_bound_f} (resp.,~\eqref{eq.Lipschitz_bound_c} for some $i \in [m]$) is not satisfied
		\State \hspace{\algorithmicindent}\hspace{\algorithmicindent} Set $L_{k,j+1} \gets \rho L_{k,j}$ (resp.,~$\gamma_{k,i,j+1} \gets \rho \gamma_{k,i,j}$)
		\State \hspace{\algorithmicindent} \textbf{else}
		\State \hspace{\algorithmicindent}\hspace{\algorithmicindent} Set $L_{k,j+1} \gets L_{k,j}$ (resp.,~$\gamma_{k,i,j+1} \gets \gamma_{k,i,j}$)
	  \EndFor
	\EndFor
  \ealgorithmic
\ealgorithm

The intuition behind the stepsize selection scheme in Algorithm~\ref{alg.sqp_adaptive} is that the stepsize is chosen to minimize an upper bound on the change in the merit function.  This upper bounding function is revealed in Lemma~\ref{lem.phi_reduction_upper} later on.  Due to the nonsmoothness of the merit function, which creates a \emph{kink} at a unit stepsize, there are three cases for the minimizer: It can occur \emph{before}, \emph{at}, or \emph{after} the kink.  An illustration of these cases is shown in Figure~\ref{fig:stepsize_illustration}.  Certain situations that lead to each of the three cases is as follows.  (There are additional situations that one may consider since the upper bounding function involves a combination of many terms, but the following are a few example situations to provide some intuition.)  If the Lipschitz constant estimates are large enough, indicating high nonlinearity of the problem functions, then the minimizer may be at a stepsize less than 1.  On the other hand, if the Lipschitz constant estimates are not too large and derivative information of the objective function suggests that the merit function improves beyond a unit stepsize, then the minimizer is at a stepsize greater than 1.  Otherwise, the minimizer occurs at a unit stepsize since at least this corresponds to a step toward linearized feasibility.

\bfigure[ht]
  \centering
  \includegraphics[width=3in,clip=true,trim=110 30 90 30]{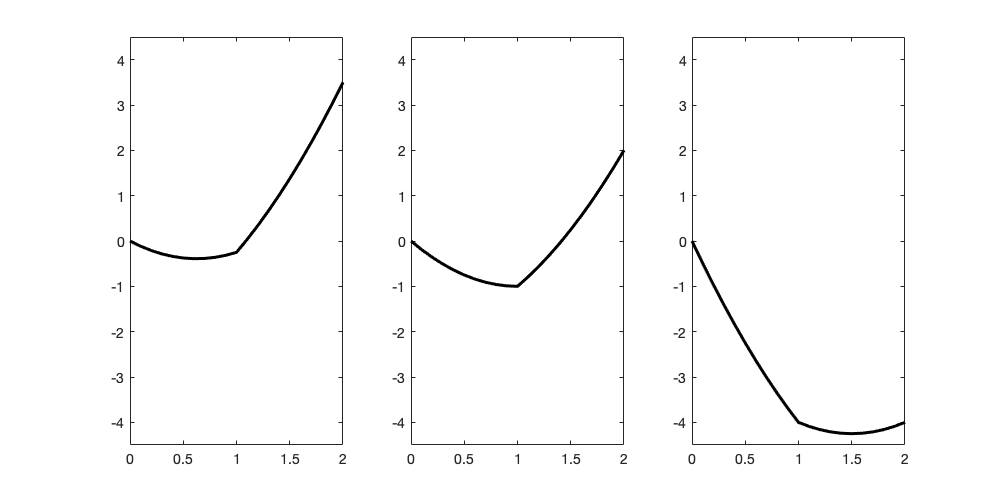}
  \caption{Illustration of three cases for an upper bounding function of the merit function (see Lemma~\ref{lem.phi_reduction_upper}) motivating the three cases in the stepsize selection scheme in Algorithm~\ref{alg.sqp_adaptive}.  Each graph shows the value of the upper bound on the change in the merit function as a function of $\alpha_k$.}
  \label{fig:stepsize_illustration}
\efigure

The second algorithm is stated as Algorithm~\ref{alg.sqp_line_search}.  In each iteration, it employs a traditional backtracking line search scheme until the reduction in the merit function is sufficiently large compared to the reduction in the model of the merit function.  This is sufficient for showing a convergence result, as shown in the next subsection.

\balgorithm[ht]
  \caption{SQP Algorithm with Backtracking Line Search}
  \label{alg.sqp_line_search}
  \balgorithmic[1]
    \Require $x_0 \in \R{n}$; $\tau_{-1} \in \R{}_{>0}$; $\epsilon \in (0,1)$; $\sigma \in (0,1)$; $\eta \in (0,1)$; $\nu \in (0,1)$; $\alpha \in \R{}_{>0}$
    \For{\textbf{all} $k \in \N{}$}
      \State Compute $(d_k,y_k)$ as the solution of \eqref{eq.system_deterministic}
      \State \textbf{if} $(x_k,y_k)$ satisfies \eqref{eq.KKT} \textbf{then return} $(x_k,y_k)$ \label{step.termination_line_search}
      \State Set $\tau_k^{trial}$ by \eqref{eq.merit_parameter_trial} and $\tau_k$ by \eqref{eq.merit_parameter_lower}
      \For{\textbf{all} $j \in \N{}$}
        \State Set $\alpha_{k,j} \gets \nu^j \alpha$
        \If{\eqref{eq.sufficient_decrease} holds}
          \State Set $\alpha_k \gets \alpha_{k,j}$ and $x_{k+1} \gets x_k + \alpha_k d_k$ and \textbf{break} (loop over $j \in \N{}$)
        \EndIf
      \EndFor
	\EndFor
  \ealgorithmic
\ealgorithm

\subsection{Convergence Analysis}

We prove in this section that, from any initial iterate, each of Algorithm~\ref{alg.sqp_adaptive} and Algorithm~\ref{alg.sqp_line_search} generates a sequence of iterates over which a first-order measure of primal-dual stationarity for \eqref{prob.deterministic} (recall \eqref{eq.KKT}) vanishes.  We assume throughout this section that both Assumptions~\ref{ass.deterministic} and \ref{ass.H} hold; for brevity, we do not remind the reader of this fact within the statement of each result.  We also remark that if an algorithm terminates finitely, then it does so with $(x_k,y_k)$ satisfying~\eqref{eq.KKT}, meaning primal-dual stationarity has been achieved.  Hence, we may assume without loss of generality in this section that neither algorithm terminates finitely, meaning that $\{x_k\}$ is infinite and $d_k \neq 0$ for all $k \in \N{}$ (recall Lemma~\ref{lem.stationary}).

In all of the results of this section, the statements are proved to hold with respect to \emph{both} Algorithms~\ref{alg.sqp_adaptive} and~\ref{alg.sqp_line_search}.  There are only a few differences in the results for the two algorithms; when a result differs, we say so explicitly.  Much of our analysis, at least prior to Lemma~\ref{lem.adaptive_red}, follows standard analysis for line-search SQP methods; see, e.g., \cite{ByrdCurtNoce08,ByrdCurtNoce10a}.  Nonetheless, we provide proofs of the results for completeness.

Our analysis uses the orthogonal decomposition of the search directions given by
\bequationNN
  d_k = u_k + v_k\ \ \text{where}\ \ u_k \in \Null(J_k)\ \ \text{and}\ \ v_k \in \Range(J_k^T)\ \ \text{for all}\ \ k \in \N{}.
\eequationNN
\emph{We emphasize that the components $u_k$ and $v_k$ do not need to be computed explicitly for any $k \in \N{}$.}  They are merely tools for our analysis.  As is common in the literature, we refer to $u_k$ as the tangential component and $v_k$ as the normal component of $d_k$.

We first show an upper bound on the normal components of the search directions.

\blemma\label{lem.bound_v}
  There exists $\kappa_v \in \R{}_{>0}$ such that, for all $k \in \N{}$, the normal component $v_k$ satisfies $\max\{\|v_k\|_2,\|v_k\|_2^2\} \leq \kappa_v \|c_k\|_2$.
\elemma
\bproof
  Let $k \in \N{}$ be arbitrary.  From $J_kd_k = J_k(u_k+v_k) = -c_k$, $u_k \in \Null(J_k)$, and $v_k \in \Range(J_k^T)$, one has $v_k = -J_k^T(J_kJ_k^T)^{-1}c_k$; hence, by Cauchy-Schwarz,
  \bequationNN
    \baligned
      \|v_k\|_2 &\leq \|J_k^T(J_kJ_k^T)^{-1}\|_2\|c_k\|_2 \\ \iff\ \ 
      \|v_k\|_2^2 &\leq (\|J_k^T(J_kJ_k^T)^{-1}\|_2\|c_k\|_2)^2 = (\|J_k^T(J_kJ_k^T)^{-1}\|_2^2\|c_k\|_2)\|c_k\|_2.
    \ealigned
  \eequationNN
  Hence, the desired conclusion follows under Assumption~\ref{ass.deterministic}.
\eproof

Our next result reveals that there exists a critical threshold between the norms of the tangential and normal components of the search directions, and in any iteration $k \in \N{}$ in which the search direction $d_k$ is dominated by the tangential component $u_k$, the curvature of $H_k$ along $d_k$ has a useful lower bound defined with $u_k$.

\blemma\label{lem.tangential_big}
  There exists $\kappa_{uv} \in \R{}_{>0}$ such that, for any $k \in \N{}$, if $\|u_k\|_2^2 \geq \kappa_{uv} \|v_k\|_2^2$, then $\thalf d_k^TH_kd_k \geq \tfrac14 \zeta \|u_k\|_2^2$, where $\zeta \in \R{}_{>0}$ is defined in Assumption~\ref{ass.H}.
\elemma
\bproof
  Assumption~\ref{ass.H} implies for any $\kappa_{uv} \in \R{}_{>0}$ that $\|u_k\|_2^2 \geq \kappa_{uv} \|v_k\|_2^2$ means
  \bequationNN
    \baligned
      \thalf d_k ^T H_k d_k &= \thalf u_k^T H_k u_k + u_k^T H_k v_k + \thalf v_k^T H_k v_k \\
      &\geq \thalf \zeta \|u_k\|_2^2 - \|u_k\|_2 \|H_k \|_2 \|v_k\|_2 - \thalf \|H_k\|_2 \|v_k\|_2^2 \\
      &\geq \(\tfrac{\zeta}{2} - \tfrac{\kappa_H}{\sqrt{\kappa_{uv}}} - \tfrac{\kappa_H}{2\kappa_{uv}}\) \|u_k\|_2^2.
    \ealigned
  \eequationNN
  Thus, under Assumption~\ref{ass.H}, the results holds for $\kappa_{uv} \in \R{}_{>0}$ with $\tfrac{\kappa_H}{\sqrt{\kappa_{uv}}} + \tfrac{\kappa_H}{2\kappa_{uv}} \leq \tfrac{\zeta}{4}$.
\eproof

For the constant $\kappa_{uv} \in \R{}_{>0}$ defined in Lemma~\ref{lem.tangential_big}, let us define
\bequationNN
  \Psi_k := \bcases \|u_k\|_2^2 + \|c_k\|_2 & \text{if $\|u_k\|_2^2 \geq \kappa_{uv} \|v_k\|_2^2$} \\ \|c_k\|_2 & \text{otherwise,} \ecases
\eequationNN
along with the related index sets (that partition $\N{}$)
\bequationNN
  \Kcal_u := \{k \in \N{} : \|u_k\|_2^2 \geq \kappa_{uv} \|v_k\|_2^2\}\ \ \text{and}\ \ \Kcal_v := \{ k \in \N{} : \|u_k\|_2^2 < \kappa_{uv} \|v_k\|_2^2\}.
\eequationNN
Our next result shows that the squared norms of the search directions and the constraint violations are bounded above by this critical quantity in all iterations.

\blemma\label{lem.Psi_1}
  There exists a constant $\kappa_\Psi \in \R{}_{>0}$ such that, for all $k \in \N{}$,
  \bequationNN
    \|d_k\|_2^2 \leq \kappa_\Psi \Psi_k\ \ \text{and}\ \ \|d_k\|_2^2 + \|c_k\|_2 \leq  (\kappa_\Psi + 1) \Psi_k.
  \eequationNN
\elemma
\bproof
  For all $k \in \Kcal_u$, it follows that
  \bequationNN
    \|d_k\|_2^2 = \|u_k\|_2^2 + \|v_k\|_2^2 \leq (1 + \kappa_{uv}^{-1}) \|u_k\|_2^2 \leq (1 + \kappa_{uv}^{-1}) (\|u_k\|_2^2 + \|c_k\|_2).
  \eequationNN
  For all $k \in \Kcal_v$, one finds from Lemma \ref{lem.bound_v} that
  \bequationNN
    \|d_k\|_2^2 = \|u_k\|_2^2 + \|v_k\|_2^2 < (\kappa_{uv} + 1) \|v_k\|_2^2 \leq (\kappa_{uv} + 1) \kappa_v \|c_k\|_2.
  \eequationNN
  Combining the results from the two cases implies the first desired result.  To establish the second result, note that the definition of $\Psi_k$ yields $\|c_k\|_2 \leq \Psi_k$ for all $k \in \N{}$.
\eproof

As revealed by our next lemma, the reduction in the model of the merit function is bounded below with respect to the same critical quantity.

\blemma\label{lem.Psi}
  There exists a constant $\kappa_q \in \R{}_{>0}$ such that, for all $k \in \N{}$,
  \bequationNN
    \Delta q(x_k,\tau_k,g_k,H_k,d_k) \geq \kappa_q \tau_k \Psi_k.
  \eequationNN
\elemma
\bproof
  Combining \eqref{eq.merit_model_reduction_lower} and Lemma~\ref{lem.tangential_big}, it follows that $\Delta q(x_k,\tau_k,g_k,H_k,d_k) \geq \tfrac14 \tau_k \zeta \|u_k\|_2^2 + \sigma \|c_k\|_1$ for $k \in \Kcal_u$. Similarly, \eqref{eq.merit_model_reduction_lower} implies that $\Delta q(x_k,\tau_k,g_k,H_k,d_k) \geq \sigma \|c_k\|_1$ for all $k \in \Kcal_v$.  Combining the two cases, $\|\cdot\|_1 \geq \|\cdot\|_2$, and the fact that $\{\tau_k\}$ is monotonically nonincreasing, the result holds for $\kappa_q := \min \{ \tfrac14 \zeta, \sigma/\tau_{-1}\} \in \R{}_{>0}$.
\eproof

Our next lemma shows an upper bound on the change in the merit function when the inner \textbf{for} loop of Algorithm~\ref{alg.sqp_adaptive} terminates with large Lipschitz constant estimates.

\blemma\label{lem.phi_reduction_upper}
  For all $k \in \N{}$, if the inner \textbf{for} loop of Algorithm~\ref{alg.sqp_adaptive} terminates since \eqref{eq.Lipschitz_bounds} holds, then with $\Gamma_k := \sum_{i\in[m]} \gamma_{k,i} \in \R{}_{>0}$ it follows that
  \bequationNN
    \phi(x_k + \alpha_k d_k,\tau_k) - \phi(x_k,\tau_k) \leq \alpha_k \tau_k g_k^Td_k + |1-\alpha_k|\|c_k\|_1 - \|c_k\|_1 + \thalf (\tau_k L_k + \Gamma_k) \alpha_k^2 \|d_k\|_2^2.
  \eequationNN
\elemma
\bproof
  For such $k \in \N{}$, it follows from \eqref{eq.Lipschitz_bounds} and Lemma~\ref{lem.c_upper} that
  \bequationNN
    \baligned
      &\ \phi(x_k + \alpha_k d_k,\tau_k) - \phi(x_k,\tau_k) \\
      =&\ \tau_k f(x_k + \alpha_k d_k) - \tau_k f(x_k) + \|c(x_k + \alpha_k d_k)\|_1 - \|c_k\|_1 \\
      \leq&\ \alpha_k \tau_k g_k^Td_k + \|c_k + \alpha_k J_k d_k\|_1 - \|c_k\|_1 + \thalf (\tau_k L_k + \Gamma_k) \alpha_k^2 \|d_k\|_2^2 \\
      =&\ \alpha_k \tau_k g_k^Td_k + |1-\alpha_k|\|c_k\|_1 - \|c_k\|_1 + \thalf (\tau_k L_k + \Gamma_k) \alpha_k^2 \|d_k\|_2^2,
    \ealigned
  \eequationNN
  as desired.
\eproof

Next, we show lower bounds for the reduction in the merit function in each iteration of each algorithm.  For concision, let us define for all $k \in \N{}$ the values
\bequationNN
  \widehat\mu_k := \tfrac{2(1-\eta)\Delta q(x_k,\tau_k,g_k,H_k,d_k) }{(\tau_k L + \sum_{i\in[m]} \gamma_i)\|d_k\|_2^2}\ \ \text{and}\ \ \widetilde\mu_k := \widehat\mu_k - \tfrac{4 \|c_k\|_1}{(\tau_k L + \sum_{i\in[m]} \gamma_{i})\|d_k\|_2^2}.
\eequationNN
For a given $k \in \N{}$, one should notice the similarity between these values and the pair $(\widehat\alpha_{k,j},\widetilde\alpha_{k,j})$ defined for all $j \in \N{}$ in Algorithm~\ref{alg.sqp_adaptive}, except that the pair $(\widehat\mu_k,\widetilde\mu_k)$ are defined with respect to $L$ and $\gamma_{i}$ for all $i \in [m]$ defined in Assumption~\ref{ass.deterministic}.

\blemma\label{lem.adaptive_red}
  For all $k \in \N{}$, the inequality \eqref{eq.sufficient_decrease} holds, where in the case of Algorithm~\ref{alg.sqp_line_search} this occurs with the stepsize satisfying $\alpha_k > \nu \min\{\widehat\mu_k,\max\{1,\widetilde\mu_k\}\} > 0$.
\elemma
\bproof
  Let $k \in \N{}$ be given.  First, consider Algorithm~\ref{alg.sqp_adaptive}.  If the inner \textbf{for} loop terminates since the stepsize yields \eqref{eq.sufficient_decrease}, then there is nothing left to prove.  Hence, we may proceed by supposing that the loop terminates since \eqref{eq.Lipschitz_bounds} holds, which we shall now proceed to show means that \eqref{eq.sufficient_decrease} holds as well.  Consider three cases, where as in Lemma~\ref{lem.phi_reduction_upper} let us define $\Gamma_k := \sum_{i\in[m]} \gamma_{k,i} \in \R{}_{>0}$.
  
  \textbf{Case 1:} Suppose that in the last iteration of the inner \textbf{for} loop one finds $\widehat\alpha_{k,j} < 1$, in which case the algorithm yields $\alpha_k = \tfrac{2(1-\eta)\Delta q(x_k,\tau_k,g_k,H_k,d_k)}{(\tau_k L_k + \Gamma_k) \|d_k\|_2^2} < 1$.  Combining this fact with Lemmas~\ref{lem.directional_derivative} and \ref{lem.phi_reduction_upper}, it follows that
  \bequationNN
    \baligned
      &\ \phi(x_k + \alpha_k d_k,\tau_k) - \phi(x_k,\tau_k) \\
      \leq&\ \alpha_k(\tau_k g_k^Td_k - \|c_k\|_1) + \thalf (\tau_k L_k + \Gamma_k) \alpha_k^2 \|d_k\|_2^2 \\
      \leq&\ - \alpha_k \Delta q(x_k,\tau_k,g_k,H_k,d_k) + \thalf (\tau_k L_k + \Gamma_k) \alpha_k^2 \|d_k\|_2^2 \\
      =&\ -\alpha_k \Delta q(x_k,\tau_k,g_k,H_k,d_k) \\
      &\quad + \thalf \alpha_k (\tau_k L_k + \Gamma_k) \(\tfrac{2(1-\eta)\Delta q(x_k,\tau_k,g_k,H_k,d_k)}{(\tau_k L_k + \Gamma_k)\|d_k\|_2^2}\)\|d_k\|_2^2 \\
      =&\ -\eta \alpha_k \Delta q(x_k,\tau_k,g_k,H_k,d_k).
    \ealigned
  \eequationNN

  \textbf{Case 2:} Suppose that in the last iteration of the inner \textbf{for} loop one finds $\widehat\alpha_{k,j} \geq 1$ and $\widetilde\alpha_{k,j} \leq 1$, in which case the algorithm yields $\alpha_k = 1$.  Combining this fact, the fact that $\widehat\alpha_{k,j} \geq 1$ in the last iteration of the loop, Lemma~\ref{lem.directional_derivative}, and Lemma~\ref{lem.phi_reduction_upper} yields the same string of relationships as in Case~1, except that since $\widehat\alpha_{k,j} \geq 1$ the first equation holds not as an equation, but as an ``$\leq$'' inequality.

  \textbf{Case 3:} Suppose that in the last iteration of the inner \textbf{for} loop one finds $\widetilde\alpha_{k,j} > 1$, in which case the algorithm yields $\alpha_k = \tfrac{2(1-\eta)\Delta q(x_k,\tau_k,g_k,H_k,d_k) - 4 \|c_k\|_1}{(\tau_k L_k + \Gamma_k) \|d_k\|_2^2} > 1$.  Combining this fact with Lemmas~\ref{lem.directional_derivative} and \ref{lem.phi_reduction_upper}, it follows that
  \bequationNN
    \baligned
      &\ \phi(x_k + \alpha_k d_k,\tau_k) - \phi(x_k,\tau_k) \\
      \leq&\ \alpha_k \tau_k g_k^Td_k + (\alpha_k-1)\|c_k\|_1 - \|c_k\|_1 + \thalf (\tau_k L_k + \Gamma_k) \alpha_k^2 \|d_k\|_2^2 \\
      =&\ \alpha_k(\tau_k g_k^Td_k - \|c_k\|_1) + 2 (\alpha_k - 1) \|c_k\|_1 + \thalf (\tau_k L_k + \Gamma_k) \alpha_k^2 \|d_k\|_2^2 \\
      \leq&\ -\alpha_k \Delta q(x_k,\tau_k,g_k,H_k,d_k) + 2\alpha_k \|c_k\|_1 + \thalf (\tau_k L_k + \Gamma_k) \alpha_k^2 \|d_k\|_2^2 \\
      =&\ -\alpha_k \Delta q(x_k,\tau_k,g_k,H_k,d_k) + 2\alpha_k \|c_k\|_1 \\
      &\quad + \thalf \alpha_k (\tau_k L_k + \Gamma_k) \(\tfrac{2(1-\eta)\Delta q(x_k,\tau_k,g_k,H_k,d_k) - 4\|c_k\|_1}{(\tau_k L_k + \Gamma_k) \|d_k\|_2^2}\) \|d_k\|_2^2 \\
      =&\ -\eta \alpha_k \Delta q(x_k,\tau_k,g_k,H_k,d_k).
    \ealigned
  \eequationNN
  Combining the three cases shows the desired result for Algorithm~\ref{alg.sqp_adaptive}.

  Now consider Algorithm~\ref{alg.sqp_line_search}.  One finds that one of three cases occurs, which mimic those for Algorithm~\ref{alg.sqp_adaptive}.  In particular, if $\widehat\mu_k < 1$, then an analysis similar to that for Case~1 above shows that for $j \in \N{}$ with $\alpha_{k,j}/\nu > \widehat\mu_k$ and $\alpha_{k,j} \leq \widehat\mu_k$, the backtracking line search will terminate by iteration $j \in \N{}$, from which it follows that $\alpha_k > \nu\widehat\mu_k$.  If $\widehat\mu_k \geq 1$ and $\widetilde\mu_k \leq 1$, or if $\widetilde\mu_k > 1$, then a similar argument combined with Case~2 or Case~3, respectively, completes the proof.
\eproof

Next, we show that the tangential components of the directions are bounded.

\blemma\label{lem.bound_u}
  The tangential component sequence $\{u_k\}$ is bounded.
\elemma
\bproof
  The first block of \eqref{eq.system_deterministic}, premultiplied by $u_k^T$, yields $u_k^TH_k(u_k + v_k) = -u_k^Tg_k$.  Hence, under Assumption \ref{ass.H}, one finds that
  \bequationNN
    \zeta \|u_k\|_2^2 \leq u_k^TH_ku_k = - g_k^Tu_k - v_k^TH_ku_k \leq (\|g_k\|_2 + \kappa_H \|v_k\|_2)\|u_k\|_2.
  \eequationNN
  Therefore, the result follows from Assumption \ref{ass.deterministic} and Lemma \ref{lem.bound_v}.
\eproof

We now show that the merit parameter sequence is bounded, and that it remains fixed at a value for all sufficiently large $k \in \N{}$.

\blemma\label{lem.tau_bound}
  There exists $k_\tau \in \N{}$ and $\tau_{\min} \in \R{}_{>0}$ such that $\tau_k = \tau_{\min}$ for $k \geq k_\tau$.
\elemma
\bproof
  Recall that $\tau_k < \tau_{k-1}$ if and only if both $g_k^Td_k + \max\{d_k^TH_kd_k,0\} > 0$ and
  \bequation\label{eq.tau_bound}
    \tau_{k-1} (g_k^Td_k + \max\{d_k^TH_kd_k,0\}) > (1 - \sigma) \|c_k\|_1.
  \eequation
  According to the first block equation of \eqref{eq.system_deterministic} (premultiplied by $u_k^T$), one has
  \bequationNN
    g_k^Td_k + \max\{d_k^TH_kd_k, 0\}
    = \bcases g_k^Tv_k + v_k^T H_k u_k + v_k^T H_k v_k & \text{if $d_k^TH_kd_k \geq 0$} \\ g_k^Tv_k - v_k^TH_ku_k - u_k^TH_ku_k & \text{otherwise.} \ecases
  \eequationNN
  The result follows from our ability to bound the left-hand side of this expression with respect to the constraint reduction.  We consider two cases.  First, if $d_k^TH_kd_k \geq 0$, then under Assumptions \ref{ass.deterministic} and \ref{ass.H} it follows with Lemmas~\ref{lem.bound_v} and \ref{lem.bound_u} and $\|\cdot\|_1 \geq \|\cdot\|_2$ that there exists a constant $\kappa_{\tau,1} \in \R{}_{>0}$ such that
  \bequationNN
    g_k^Tv_k + v_k^T H_k u_k + v_k^T H_k v_k \leq (\|g_k\|_2 + \kappa_H \|u_k\|_2) \|v_k\|_2 + \kappa_H \|v_k\|_2^2 \leq \kappa_{\tau,1} \|c_k\|_1.
  \eequationNN
  Second, if $d_k^TH_kd_k < 0$, then under Assumptions \ref{ass.deterministic} and \ref{ass.H} it follows from Lemmas \ref{lem.bound_v} and \ref{lem.bound_u} and $\|\cdot\|_1 \geq \|\cdot\|_2$ that there exists a constant $\kappa_{\tau,2} \in \R{}_{>0}$ such that
  \bequationNN
    g_k^Tv_k - v_k^T H_k u_k - u_k^T H_k u_k \leq (\|g_k\|_2 + \kappa_H \|u_k\|_2) \|v_k\|_2 \leq \kappa_{\tau,2} \|c_k\|_1.
  \eequationNN
  Together, these imply $g_k^Td_k + \max\{d_k^TH_kd_k,0\} \leq \max\{\kappa_{\tau,1},\kappa_{\tau,2}\} \|c_k\|_1$, from which it follows that in order to have both $g_k^Td_k + \max\{d_k^TH_kd_k,0\} > 0$ and \eqref{eq.tau_bound}, one must have $\tau_{k-1} > (1-\sigma)/\max\{\kappa_{\tau,1},\kappa_{\tau,2}\}$.  Therefore, if this inequality is not satisfied for $k = k_\tau$ for some $k_\tau \in \R{}$, then it remains unsatisfied for all $k \geq k_\tau$.  This, along with the fact that whenever Algorithm \ref{alg.sqp_adaptive} or \ref{alg.sqp_line_search} decreases the merit parameter it does so by at least a constant factor, proves the result.
\eproof

We now prove that there is a positive lower bound for the stepsizes.

\blemma\label{lem.alpha_lower}
  There exists $\alpha_{\min} \in \R{}_{>0}$ such that $\alpha_k \geq \alpha_{\min}$ for all $k \in \N{}$.
\elemma
\bproof
  Let $k \in \N{}$ be given.  With respect to Algorithm~\ref{alg.sqp_adaptive}, one has that $\alpha_k \geq 1$ unless the inner \textbf{for} loop terminates in iteration $j \in \N{}$ with $\widehat\alpha_{k,j} < 1$.  In such cases, it follows from monotonicity of $\{\tau_k\}$ and Lemmas~\ref{lem.well_defined}, \ref{lem.tau_bound}, \ref{lem.Psi_1}, and \ref{lem.Psi} that
  \bequationNN
    \alpha_k = \tfrac{2(1-\eta) \Delta q(x_k,\tau_k,g_k,H_k,d_k)}{(\tau_k L_{k,j} + \sum_{i\in[m]} \gamma_{k,i,j})\|d_k\|_2^2} \geq \tfrac{2(1-\eta)\kappa_q \tau_{\min}}{(\tau_{-1} L_{\max} + \sum_{i\in[m]} \gamma_{\max,i}) \kappa_\Psi} > 0.
  \eequationNN
  Similarly, for Algorithm~\ref{alg.sqp_line_search}, Lemma~\ref{lem.adaptive_red} implies $\alpha_k \geq 1$ unless $\widehat\mu_k < 1$.  In such cases, it follows from monotonicity of $\{\tau_k\}$ and Lemmas~\ref{lem.well_defined}, \ref{lem.tau_bound}, \ref{lem.Psi_1}, and \ref{lem.Psi} that
  \bequationNN
    \alpha_k > \tfrac{2\nu (1-\eta) \Delta q(x_k,\tau_k,g_k,H_k,d_k)}{(\tau_k L + \sum_{i\in[m]} \gamma_i)\|d_k\|_2^2} \geq \tfrac{2 \nu (1-\eta)\kappa_q \tau_{\min}}{(\tau_{-1} L + \sum_{i\in[m]} \gamma_{i}) \kappa_\Psi} > 0.
  \eequationNN
  Overall, a positive lower bound has been proved for both algorithms.
\eproof

We now present our main convergence theorem for Algorithms~\ref{alg.sqp_adaptive} and \ref{alg.sqp_line_search}.

\btheorem\label{th.deterministic}
  Algorithms~\ref{alg.sqp_adaptive} and \ref{alg.sqp_line_search} yield
  \bequationNN
    \lim_{k\to\infty} \|d_k\|_2 = 0,\ \ \lim_{k\to\infty} \|c_k\|_2 = 0,\ \ \text{and}\ \ \lim_{k\to\infty} \|g_k + J_k^Ty_k\|_2 = 0.
  \eequationNN    
\etheorem
\bproof
  For all $k \in \N{}$, it follows from Lemmas~\ref{lem.Psi}, \ref{lem.adaptive_red}, and \ref{lem.alpha_lower} that
  \bequationNN
    \phi(x_k,\tau_k) - \phi(x_{k+1},\tau_k) \geq \eta \alpha_k \Delta q(x_k,\tau_k,g_k,H_k,d_k) \geq \eta \alpha_{\min} \kappa_q \tau_{\min} \Psi_k.
  \eequationNN
  Combining this with Lemmas~\ref{lem.Psi_1} and \ref{lem.tau_bound} shows for $k \in \N{}$ with $k > k_\tau$ that
  \bequationNN
    \baligned
      &\ \phi(x_{k_\tau},\tau_{\min}) - \phi(x_k,\tau_{\min}) \\
      =&\ \sum_{j = k_\tau}^{k-1} (\phi(x_j,\tau_{\min}) - \phi(x_{j+1},\tau_{\min})) \\
      \geq&\ \eta \alpha_{\min} \kappa_q \tau_{\min} \sum_{j = k_\tau}^{k-1} \Psi_j \geq \tfrac{\eta \alpha_{\min} \kappa_q \tau_{\min}} {\kappa_\Psi + 1} \sum_{j = k_\tau}^{k-1} (\|d_j\|_2^2 + \|c_j\|_2).
    \ealigned
  \eequationNN
  Since, under Assumption~\ref{ass.deterministic}, $\phi(\cdot,\tau_{\min})$ is bounded below over the iterates, the above implies the first two desired limits.  Note now that \eqref{eq.system_deterministic} implies
  \bequation\label{eq.gJy}
    \|g_k + J_k^Ty_k\|_2 = \|H_kd_k\|_2 \leq \|H_k\|_2\|d_k\|_2 \leq \kappa_H \|d_k\|_2.
  \eequation
  Hence, by Assumption~\ref{ass.H} and $\{d_k\} \to 0$, the result follows.
\eproof

\section{Stochastic Setting}\label{sec.stochastic}

Now consider the optimization problem
\bequation\label{prob.f_nonlinear_stochastic}
  \min_{x\in\R{n}}\ f(x)\ \st\ c(x) = 0,\ \ \text{with}\ \ f(x) = \E[F(x,\omega)],
\eequation
where $f : \R{n} \to \R{}$, $c : \R{n} \to \R{m}$, $\omega$ is a random variable with associated probability space $(\Omega,\Fcal,P)$, $F : \R{n} \times \Omega \to \R{}$, and $\E[\cdot]$ represents expectation taken with respect to~$P$.  We presume that one has access to values of the constraint function and its derivatives, but that it is intractable to evaluate the objective and/or its derivatives.  That said, we presume that at a given iterate $x_k$, one can evaluate a stochastic gradient estimate $\gbar_k \in \R{n}$ satisfying the following assumption.

\bassumption\label{ass.g}
  For all $k \in \N{}$, the stochastic gradient estimate $\gbar_k \in \R{n}$ is an unbiased estimator of the gradient of $f$ at $x_k$, i.e.,
  \bequationNN
    \E_k[\gbar_k] = g_k,
  \eequationNN
  where $\E_k[\cdot]$ denotes expectation taken with respect to the distribution of $\omega$ conditioned on the event that the algorithm has reached $x_k \in \R{n}$ in iteration $k \in \N{}$.  In addition, there exists a constant $M \in \R{}_{>0}$ such that, for all $k \in \N{}$, one has
  \bequationNN
    \E_k[\|\gbar_k - g_k\|_2^2] \leq M.
  \eequationNN
\eassumption

\subsection{Algorithm}\label{sec.algorithms_stochastic}

Similar to the deterministic setting, in order to solve \eqref{prob.f_nonlinear_stochastic}, we consider a stochastic algorithm that computes a search direction $\dbar_k \in \R{n}$ and Lagrange multiplier vector $\ybar_k \in \R{m}$ in iteration $k \in \N{}$ by solving the linear system
\bequation\label{eq.system_stochastic}
  \bbmatrix H_k & J_k^T \\ J_k & 0 \ebmatrix \bbmatrix \dbar_k \\ \ybar_k \ebmatrix = - \bbmatrix \gbar_k \\ c_k \ebmatrix,
\eequation
where $\{H_k\}$ satisfies Assumption~\ref{ass.H}.  Generally, we use a ``bar'' over a quantity whose value in iteration $k \in \N{}$ depends on $\gbar_k$.  Hence, as they are independent of $\gbar_k$ conditioned on the event that the algorithm reaches $x_k$ as its $k$th iterate, we write the constraint value, constraint Jacobian, and $(1,1)$-block matrix as $c_k$, $J_k$, and $H_k$, respectively, but we write the solution of \eqref{eq.system_stochastic} as $(\dbar_k,\ybar_k)$ due to its dependence on~$\gbar_k$.

The algorithm that we propose is stated as Algorithm~\ref{alg.sqp_stochastic}.  Paralleling Algorithm~\ref{alg.sqp_adaptive}, the merit parameter is set based on the computation of a trial value
\bequation\label{eq.merit_parameter_trial_stochastic}
  \bar\tau_k^{trial} \gets \bcases \infty & \text{if $\gbar_k^T\dbar_k + \max\{\dbar_k^TH_k\dbar_k,0\} \leq 0$} \\ \tfrac{(1 - \sigma)\|c_k\|_1}{\gbar_k^T\dbar_k + \max\{\dbar_k^TH_k\dbar_k,0\}} & \text{otherwise,} \ecases
\eequation
followed by the rule
\bequation\label{eq.merit_parameter_lower_stochastic}
  \bar\tau_k \gets \bcases \bar\tau_{k-1} & \text{if $\bar\tau_{k-1} \leq \bar\tau_k^{trial}$} \\ (1-\epsilon) \bar\tau_k^{trial} & \text{otherwise,} \ecases
\eequation
which ensures $\bar\tau_k \leq \bar\tau_k^{trial}$ and, similarly as for our deterministic algorithm (see \eqref{eq.merit_model_reduction_lower}),
\bequation\label{eq.merit_model_reduction_lower_stochastic}
  \Delta q(x_k,\bar\tau_k,\gbar_k,H_k,\dbar_k) \geq \thalf \bar\tau_k \max\{\dbar_k^TH_k\dbar_k,0\} + \sigma \|c_k\|_1.
\eequation

A unique feature of our algorithm for this stochastic setting is that it adaptively estimates a lower bound for the ratio between the reduction in the model of the merit function and the merit parameter times the squared norm of a search direction.  This is used to determine an interval into which the stepsize will be projected; control of this parameter is paramount to ensure convergence in expectation.  We set
\bequation\label{eq.ratio_trial}
  \bar\xi_k^{trial} \gets \tfrac{\Delta q(x_k,\bar\tau_k,\gbar_k,H_k,\dbar_k)}{\bar\tau_k \|\dbar_k\|_2^2},
\eequation
then apply the rule (which ensures $\bar\xi_k \leq \bar\xi_k^{trial}$ for all $k \in \N{}$)
\bequation\label{eq.ratio}
  \bar\xi_k \gets \bcases \bar\xi_{k-1} & \text{if $\bar\xi_{k-1} \leq \bar\xi_k^{trial}$} \\ (1-\epsilon) \bar\xi_k^{trial} & \text{otherwise.} \ecases
\eequation
It will be shown in our analysis that $\{\bar\xi_k\}$ is bounded away from zero \emph{deterministically}.

For generality, Algorithm~\ref{alg.sqp_stochastic} is stated with Lipschitz constant estimates $\{L_k\}$ and~$\{\Gamma_k\}$ given as inputs (with the idea that $\Gamma_k := \sum_{i\in[m]} \gamma_{k,i}$ for all $k \in \N{}$).  Our analysis in the next subsection presumes that Lipschitz constants are known, although in practice these can be estimated using standard techniques (see, e.g., \cite{CurtRobi19}) in an attempt to ensure that the same convergence results hold as for the case when the constants are known.  The sequence $\{\beta_k\}$ is introduced to control the stepsizes.  As in standard analysis for stochastic (sub)gradient-type methods, our analysis in the next subsection considers the case when $\{\beta_k\}$ is constant asymptotically, and when it diminishes at an appropriate rate to ensure convergence in expectation.

\begin{algorithm}[ht]
  \caption{Stochastic SQP Algorithm}
  \label{alg.sqp_stochastic}
  \begin{algorithmic}[1]
    \Require $x_0 \in \R{n}$; $\bar\tau_{-1} \in \R{}_{>0}$; $\epsilon \in (0,1)$; $\sigma \in (0,1)$; $\bar\xi_{-1} \in \R{}_{>0}$; $\{\beta_k\} \subset (0,1]$; $\theta \in \R{}_{\geq0}$; $\{L_k\} \subset \R{}_{>0}$; $\{\Gamma_k\} \subset \R{}_{>0}$
    \For{\textbf{all} $k \in \N{}$}
	  \State Compute $(\dbar_k,\ybar_k)$ as the solution of \eqref{eq.system_stochastic}
	  \State \textbf{if} $\dbar_k = 0$ \textbf{then continue} (to iteration $k+1$)
	  \State Set $\bar\tau_k^{trial}$ by \eqref{eq.merit_parameter_trial_stochastic} and $\bar\tau_k$ by \eqref{eq.merit_parameter_lower_stochastic} \label{step.tau_stochastic}
	  \State Set $\bar\xi_k^{trial}$ by \eqref{eq.ratio_trial} and $\bar\xi_k$ by \eqref{eq.ratio} \label{step.xi}
      \State Set
      \bequationNN
        \baligned
          \bar{\widehat\alpha}_{k,\text{init}} &\gets \tfrac{\beta_k\Delta q(x_k,\bar\tau_k,\gbar_k,H_k,\dbar_k)}{(\bar\tau_k L_k + \Gamma_k)\|\dbar_k\|_2^2}\ \ \text{and} \\ \bar{\widetilde\alpha}_{k,\text{init}} &\gets \bar{\widehat\alpha}_{k,\text{init}} - \tfrac{4\|c_k\|_1}{(\bar\tau_k L_k +  \Gamma_k)\|\dbar_k\|_2^2}
        \ealigned
      \eequationNN
      \State Set $\bar{\widehat\alpha}_k \gets \proj_k(\bar{\widehat\alpha}_{k,\text{init}})$ and $\bar{\widetilde\alpha}_k \gets \proj_k(\bar{\widetilde\alpha}_{k,\text{init}})$ where \label{step.alpha_projection_stochastic}
      \bequationNN
        \proj_k(\cdot) \equiv \proj\( \cdot\ \bigg| \left[\tfrac{\beta_k \bar\xi_k \bar\tau_k}{\bar\tau_k L_k + \Gamma_k}, \tfrac{\beta_k \bar\xi_k \bar\tau_k}{\bar\tau_k L_k + \Gamma_k} + \theta\beta_k^2\right]\)
      \eequationNN
      \State Set \label{step.alpha_stochastic}
	  \bequationNN
	    \bar{\alpha}_k \gets
		  \bcases
		    \bar{\widehat\alpha}_k & \text{if $\bar{\widehat\alpha}_k < 1$} \\
		    1 & \text{if $\bar{\widetilde\alpha}_k \leq 1 \leq \bar{\widehat\alpha}_k$} \\
		    \bar{\widetilde\alpha}_k & \text{if $\bar{\widetilde\alpha}_k > 1$}
		  \ecases
		\eequationNN
      \State Set $x_{k+1} \gets x_k + \bar\alpha_k \dbar_k$
    \EndFor
  \end{algorithmic}
\end{algorithm}

\subsection{Convergence Analysis}\label{sec.analysis_stochastic}

In this section, we prove that Algorithm~\ref{alg.sqp_stochastic} has convergence properties that match those from the deterministic setting in expectation, with some caveats that we explain and justify.  Our algorithm uses only the stochastic gradient estimates~$\{\gbar_k\}$, computes $\{(\dbar_k,\ybar_k)\}$ by \eqref{eq.system_stochastic}, sets merit parameter-related sequences $\{\bar\tau_k\}$ and~$\{\bar\tau_k^{trial}\}$, and also sets steplength-related sequences $\{\bar\xi_k\}$ and~$\{\bar\xi_k^{trial}\}$, but our analysis also references the gradients $\{g_k\}$ corresponding to $\{x_k\}$ as well as the corresponding sequence of solutions of \eqref{eq.system_deterministic}, namely, $\{(d_k,y_k)\}$, and trial merit parameter values $\{\tau_k^{trial}\}$.  In other words, for all $k \in \N{}$, conditioned on the event that the algorithm reaches $x_k$, we define $(d_k,y_k)$ and $\tau_k^{trial}$ as they would be computed if the algorithm reached $x_k$ as the $k$th iterate in Algorithm~\ref{alg.sqp_adaptive}.  Throughout this section, we assume that Assumptions~\ref{ass.deterministic}, \ref{ass.H}, and \ref{ass.g} hold---where $\{H_k\}$ is a deterministic sequence chosen independently from $\{\gbar_k\}$---and for the sake of brevity we do not state this fact within each result.  In this section, we assume that Lipschitz constants for the objective and constraints, in particular, $L$ and $\Gamma := \sum_{i\in[m]} \gamma_{i}$, are known.

\begin{remark}
  Our analysis makes Assumption~\ref{ass.deterministic}, which means that it assumes that the iterates remain in an open convex set over which the objective and constraint function and derivative values are bounded.  This is admittedly not ideal in a stochastic setting.  For example, in the case of applying a stochastic gradient method (SG) in an unconstrained stochastic setting, it is not ideal to assume that the gradients at the iterates remain bounded in norm, since---as SG is not a descent method---it is unreasonable to assume that the iterates remain in a sublevel set of the objective function.  However, we believe this assumption is more reasonable in a constrained setting, since the iterates are being driven to the \emph{deterministic} feasible region.  Further, we claim that Assumption~\ref{ass.deterministic} could be loosened if our algorithm were to choose a predetermined stepsize sequence, rather than one that mimicks the stepsize scheme from Algorithm~\ref{alg.sqp_adaptive}.  We discuss this issue further in \S\ref{sec.conclusion}.
\end{remark}

Since if $\dbar_k = 0$, then the algorithm simply skips to iteration $k+1$, we may assume without loss of generality in our analysis that $\dbar_k \neq 0$ for all $k \in \N{}$.

As in the deterministic setting, our analysis makes use of the orthogonal decomposition of the (stochastic) search directions given by
\bequationNN
  \dbar_k = \bar{u}_k + v_k\ \ \text{where}\ \ \bar{u}_k \in \Null(J_k)\ \ \text{and}\ \ v_k \in \Range(J_k^T)\ \ \text{for all}\ \ k \in \N{}.
\eequationNN
Let us emphasize that, conditioned on the event that the algorithm reaches $x_k$ as its $k$th iterate, the normal component is \emph{deterministic}, depending only on the constraint value $c_k$ and Jacobian $J_k$; hence, we write $v_k$ rather than~$\vbar_k$ in the expression above.  For all $k \in \N{}$, let $Z_k$ be an orthogonal basis for the null space of~$J_k$, which under Assumption~\ref{ass.deterministic} is a matrix in $\R{n \times (n-m)}$.  It follows that, for all $k \in \N{}$,
\bequationNN
  \ubar_k = Z_k\wbar_k\ \ \text{and}\ \ u_k = Z_kw_k\ \ \text{for some}\ \ (\wbar_k,w_k) \in \R{n-m} \times \R{n-m}.
\eequationNN
Under Assumption~\ref{ass.H}, the reduced Hessian satisfies $Z_k^TH_kZ_k \succeq \zeta I$.

For our first lemma, we carry over properties of algorithmic quantities that hold in the same manner as in the deterministic case, conditioned on the event that the algorithm has reached $x_k$ as the $k$th iterate.  As in our analysis in the deterministic setting, for the constant $\kappa_{uv} \in \R{}_{>0}$ defined in the lemma, we define
\bequationNN
  \overline\Psi_k := \bcases \|\ubar_k\|_2^2 + \|c_k\|_2 & \text{if $\|\ubar_k\|_2^2 \geq \kappa_{uv} \|v_k\|_2^2$} \\ \|c_k\|_2 & \text{otherwise.} \ecases
\eequationNN

\blemma\label{lem.deterministic_to_stochastic}
  For all $k \in \N{}$, \eqref{eq.system_stochastic} has a unique solution.  In addition, for the same constants $(\kappa_v,\kappa_{uv},\kappa_\Psi,\kappa_q) \in \R{}_{>0} \times \R{}_{>0} \times \R{}_{>0} \times \R{}_{>0}$ that appear in Lemmas~\ref{lem.bound_v}, \ref{lem.tangential_big}, \ref{lem.Psi_1}, and \ref{lem.Psi}, the following statements hold true for all $k \in \N{}$.
  \benumerate
    \item[(a)] The normal component satisfies $\max\{\|v_k\|_2,\|v_k\|_2^2\} \leq \kappa_v\|c_k\|_2$.
    \item[(b)] If $\|\ubar_k\|_2^2 \geq \kappa_{uv}\|v_k\|_2^2$, then $\thalf \dbar_k^TH_k\dbar_k \geq \tfrac14 \zeta \|\ubar_k\|_2^2$.
    \item[(c)] The search direction satisfies $\|\dbar_k\|_2^2 \leq \kappa_\Psi \overline\Psi_k$ and $\|\dbar_k\|_2^2 + \|c_k\|_2 \leq (\kappa_\Psi + 1) \overline\Psi_k$.
    \item[(d)] The model reduction satisfies $\Delta q(x_k,\bar\tau_k,\gbar_k,H_k,\dbar_k) \geq \kappa_q \bar\tau_k \overline\Psi_k$.
  \eenumerate
  Finally, for all $k \in \N{}$, it follows that
  \bequationNN
    \phi(x_k + \bar\alpha_k \dbar_k,\bar\tau_k) - \phi(x_k,\bar\tau_k) \leq \bar\alpha_k \bar\tau_k g_k^T \dbar_k + |1-\bar\alpha_k| \|c_k\|_1 - \|c_k\|_1 + \thalf (\bar\tau_k L + \Gamma) \bar\alpha_k^2 \|\dbar_k\|_2^2.
  \eequationNN
\elemma
\bproof
  That \eqref{eq.system_stochastic} has a unique solution for all $k \in \N{}$ follows for the same reason that Lemma~\ref{lem.nonsingular} holds.  The proofs of parts (a)--(d) follow in the same manner as the proofs of Lemmas~\ref{lem.bound_v}, \ref{lem.tangential_big}, \ref{lem.Psi_1}, and \ref{lem.Psi}, respectively, with the stochastic quantities $\{\gbar_k,\dbar_k,\ubar_k,\bar\tau_k\}$ in place of the deterministic quantities $\{g_k,d_k,u_k,\tau_k\}$, where it is important to recognize that the conclusions follow with the \emph{same constants}, namely, $(\kappa_v,\kappa_{uv},\kappa_\Psi,\kappa_q)$, as in the deterministic setting.  The proof of the last conclusion follows in the same manner as that of Lemma~\ref{lem.phi_reduction_upper}.
\eproof

In the next lemma, we prove that the sequence $\{\bar\xi_k\}$ is bounded deterministically.

\blemma\label{lem.xi}
  In any run of the algorithm, there exists $\kbar_\xi \in \N{}$ and $\bar\xi_{\min} \in \R{}_{>0}$ such that $\bar\xi_k = \bar\xi_{\min}$ for all $k \geq \kbar_\xi$, where $\bar\xi_{\min} \in [\xi_{\min},\bar\xi_{-1}]$ with $\xi_{\min} := \epsilon\kappa_q/\kappa_\Psi$.
\elemma
\bproof
  If Line~\ref{step.xi} of the algorithm ever sets $\bar\xi_k < \bar\xi_{k-1}$, then it ensures that $\bar\xi_k \leq \epsilon \bar\xi_{k-1}$.  This means that $\{\bar\xi_k\}$ is constant for sufficiently large $k$ or it vanishes.  On the other hand, by Lemma~\ref{lem.deterministic_to_stochastic}(c) and (d), it follows that
  \bequationNN
    \tfrac{\Delta q(x_k,\bar\tau_k,\gbar_k,H_k,\dbar_k)}{\bar\tau_k \|\dbar_k\|_2^2} \geq \tfrac{\kappa_q \bar\tau_k \overline\Psi_k}{\kappa_\Psi \bar\tau_k \overline\Psi_k} = \tfrac{\kappa_q}{\kappa_\Psi},
  \eequationNN
  meaning that Line~\ref{step.xi} will never set $\bar\xi_k$ less than $\epsilon\kappa_q/\kappa_\Psi$ for any $k \in \N{}$.  Therefore, $\{\bar\xi_k\}$ is constant for sufficiently large $k$ in the manner stated.
\eproof

Next, we present the following obvious, but important consequence of our stepsize selection scheme.  In particular, the result shows that, even though the algorithm sets the stepsize adaptively, the difference between the largest and smallest possible stepsizes in a given iteration is $\Ocal(\beta_k^2)$, so this difference is controlled by the algorithm.

\blemma\label{lem.beta_control}
  For any $k \in \N{}$, the stepsize satisfies
  \bequationNN
    \bar\alpha_k \in [\bar\alpha_{k,\min},\bar\alpha_{k,\max}] := \left[\tfrac{\beta_k \bar\xi_k \bar\tau_k }{\bar\tau_k L + \Gamma}, \tfrac{\beta_k \bar\xi_k \bar\tau_k}{\bar\tau_k L + \Gamma} + \theta\beta_k^2\right],
  \eequationNN
  which is an interval with length $\bar\alpha_{k,\max} - \bar\alpha_{k,\min} = \theta \beta_k^2$.
\elemma
\bproof
  The proof follows directly from the projections of $\bar{\widehat\alpha}_k$ and $\bar{\widetilde\alpha}_k$ in  Line~\ref{step.alpha_projection_stochastic} and the formula for the stepsize $\bar\alpha_k$ in Line~\ref{step.alpha_stochastic}.
\eproof

Our next result is a cornerstone of our analysis.  It builds on the last conclusion in Lemma~\ref{lem.deterministic_to_stochastic} to specify a useful upper bound for the merit function value after a step.  Central to the proof is our specific stepsize selection strategy.

\blemma\label{lem.key_decrease}
  Suppose that $\{\beta_k\}$ is chosen such that $\beta_k \bar\xi_k \bar\tau_k/(\bar\tau_k L + \Gamma) \in (0,1]$ for all $k \in \N{}$.  Then, for all $k \in \N{}$, it follows that
  \bequationNN
    \baligned
      &\ \phi(x_k + \bar\alpha_k \dbar_k, \bar\tau_k) - \phi(x_k, \bar\tau_k) \\
      \leq&\ -\bar\alpha_k \Delta q(x_k,\bar\tau_k,g_k,H_k,d_k) + \thalf \bar\alpha_k \beta_k  \Delta q(x_k,\bar\tau_k,\gbar_k,H_k,\dbar_k) + \bar\alpha_k \bar\tau_k g_k^T (\dbar_k - d_k).
    \ealigned
  \eequationNN
\elemma
\bproof
  Let $k \in \N{}$ be arbitrary.  We consider three cases, with a few subcases, depending on how the stepsize is set in Lines~\ref{step.alpha_projection_stochastic} and~\ref{step.alpha_stochastic} of the algorithm.

  \textbf{Case 1:} Suppose in Line~\ref{step.alpha_stochastic} that $\bar{\widehat\alpha}_k < 1$, meaning that $\bar\alpha_k \gets \bar{\widehat\alpha}_k$.  From Lemma~\ref{lem.deterministic_to_stochastic} and Lemma~\ref{lem.directional_derivative}, it follows that
  \bequationNN
    \baligned
      &\ \phi(x_k + \bar\alpha_k \dbar_k, \bar\tau_k) - \phi(x_k, \bar\tau_k) \\
      \leq&\ \bar\alpha_k (\bar\tau_k g_k^T \dbar_k - \|c_k\|_1) + \thalf (\bar\tau_k L + \Gamma) \bar\alpha_k^2 \|\dbar_k\|_2^2 \\
      =&\ \bar\alpha_k (\bar\tau_k g_k^T d_k - \|c_k\|_1) + \thalf (\bar\tau_k L + \Gamma) \bar\alpha_k^2 \|\dbar_k\|_2^2 + \bar\alpha_k \bar\tau_k g_k^T (\dbar_k - d_k) \\
      \leq&\ -\bar\alpha_k \Delta q(x_k,\bar\tau_k,g_k,H_k,d_k) + \thalf (\bar\tau_k L + \Gamma) \bar\alpha_k^2 \|\dbar_k\|_2^2 + \bar\alpha_k \bar\tau_k g_k^T (\dbar_k - d_k).
    \ealigned
  \eequationNN
  Using this inequality, let us now consider two subcases.  (For all $k \in \N{}$, since \eqref{eq.ratio} ensures $\bar\xi_k \leq \bar\xi_k^{trial} = \tfrac{\Delta q(x_k,\bar\tau_k,\gbar_k,H_k,\dbar_k)}{\bar\tau_k\|\dbar_k\|_2^2}$, it follows that $\tfrac{\beta_k\Delta q(x_k,\bar\nu_k,\gbar_k,H_k,\dbar_k)}{(\bar\tau_k L + \Gamma)\|\dbar_k\|_2^2} \geq \tfrac{\beta_k\bar\xi_k\bar\tau_k}{\bar\tau_k L + \Gamma}$.)
  
  \noindent
  \textbf{Case 1a:} If $\bar\alpha_k = \tfrac{\beta_k\Delta q(x_k,\bar\tau_k,\gbar_k,H_k,\dbar_k)}{(\bar\tau_k L + \Gamma)\|\dbar_k\|_2^2}$, then
  \bequationNN
    \baligned
      &\ \phi(x_k + \bar\alpha_k \dbar_k, \bar\tau_k) - \phi(x_k, \bar\tau_k) \\
      \leq&\ -\bar\alpha_k \Delta q(x_k,\bar\tau_k,g_k,H_k,d_k) \\
      &\quad + \thalf \bar\alpha_k (\bar\tau_k L + \Gamma) \(\tfrac{\beta_k \Delta q(x_k,\bar\tau_k,\gbar_k,H_k,\dbar_k)}{(\bar\tau_k L + \Gamma) \|\dbar_k\|_2^2}\) \|\dbar_k\|_2^2 + \bar\alpha_k \bar\tau_k g_k^T (\dbar_k - d_k) \\
      =&\ -\bar\alpha_k \Delta q(x_k,\bar\tau_k,g_k,H_k,d_k) + \thalf \bar\alpha_k \beta_k \Delta q(x_k,\bar\tau_k,\gbar_k,H_k,\dbar_k) + \bar\alpha_k \bar\tau_k g_k^T (\dbar_k - d_k).
    \ealigned
  \eequationNN
  
  \noindent
  \textbf{Case 1b:} If $\bar\alpha_k = \tfrac{\beta_k \bar\xi_k \bar\tau_k}{\bar\tau_k L + \Gamma} + \theta\beta_k^2 \leq \tfrac{\beta_k \Delta q(x_k,\bar\tau_k,\gbar_k,H_k,\dbar_k)}{(\bar\tau_k L + \Gamma) \|\dbar_k\|_2^2}$, then
  \bequationNN
    \baligned
      &\ \phi(x_k + \bar\alpha_k \dbar_k, \bar\tau_k) - \phi(x_k, \bar\tau_k) \\
      \leq&\ -\bar\alpha_k \Delta q(x_k,\bar\tau_k,g_k,H_k,d_k) \\
      &\quad + \thalf \bar\alpha_k (\bar\tau_k L + \Gamma) \(\tfrac{\beta_k \bar\xi_k \bar\tau_k}{\bar\tau_k L + \Gamma} + \theta\beta_k^2\) \|\dbar_k\|_2^2 + \bar\alpha_k \bar\tau_k g_k^T (\dbar_k - d_k) \\
      \leq&\ - \bar\alpha_k \Delta q(x_k,\bar\tau_k,g_k,H_k,d_k) + \thalf \bar\alpha_k \beta_k \Delta q(x_k,\bar\tau_k,\gbar_k,H_k,\dbar_k) + \bar\alpha_k \bar\tau_k g_k^T (\dbar_k - d_k).
    \ealigned
  \eequationNN
  
  
  \textbf{Case 2:} Suppose in Line~\ref{step.alpha_stochastic} that $\bar{\widetilde\alpha}_k \leq 1 \leq \bar{\widehat\alpha}_k$, meaning that $\bar\alpha_k \gets 1$.  From Lemma~\ref{lem.deterministic_to_stochastic}, Lemma~\ref{lem.directional_derivative}, and since $\tfrac{\beta_k \Delta q(x_k,\bar\tau_k,\gbar_k,H_k,\dbar_k)}{(\bar\tau_k L + \Gamma) \|\dbar_k\|_2^2} \geq 1 = \bar\alpha_k$, it follows that
  \bequationNN
    \baligned
      &\ \phi(x_k + \bar\alpha_k \dbar_k, \bar\tau_k) - \phi(x_k, \bar\tau_k) \\
      \leq&\ \bar\alpha_k (\bar\tau_k g_k^T \dbar_k - \|c_k\|_1) + \thalf (\bar\tau_k L + \Gamma) \bar\alpha_k^2 \|\dbar_k\|_2^2 \\
      =&\ \bar\alpha_k (\bar\tau_k g_k^T d_k - \|c_k\|_1) + \thalf (\bar\tau_k L + \Gamma) \bar\alpha_k^2 \|\dbar_k\|_2^2 + \bar\alpha_k \bar\tau_k g_k^T (\dbar_k - d_k) \\
      \leq&\ -\bar\alpha_k \Delta q(x_k,\bar\tau_k,g_k,H_k,d_k) + \thalf (\bar\tau_k L + \Gamma) \bar\alpha_k^2 \|\dbar_k\|_2^2 + \bar\alpha_k \bar\tau_k g_k^T (\dbar_k - d_k) \\
      \leq&\ -\bar\alpha_k \Delta q(x_k,\bar\tau_k,g_k,H_k,d_k) + \thalf \bar\alpha_k \beta_k \Delta q(x_k,\bar\tau_k,\gbar_k,H_k,\dbar_k) + \bar\alpha_k \bar\tau_k g_k^T (\dbar_k - d_k).
    \ealigned
  \eequationNN
  
  \textbf{Case 3:} Suppose in Line~\ref{step.alpha_stochastic} that $\bar{\widetilde\alpha}_k > 1$, meaning that $\bar\alpha_k \gets \bar{\widetilde\alpha}_k$.  From Lemma~\ref{lem.deterministic_to_stochastic} and Lemma~\ref{lem.directional_derivative}, it follows that
  \bequationNN
    \baligned
      &\ \phi(x_k + \bar\alpha_k \dbar_k, \bar\tau_k) - \phi(x_k, \bar\tau_k) \\
      \leq&\ \bar\alpha_k \bar\tau_k g_k^T \dbar_k + (\bar\alpha_k - 1)\|c_k\|_1 - \|c_k\|_1 + \thalf (\bar\tau_k L + \Gamma) \bar\alpha_k^2 \|\dbar_k\|_2^2 \\
      =&\ \bar\alpha_k (\bar\tau_k g_k^T \dbar_k - \|c_k\|_1) + 2(\bar\alpha_k - 1)\|c_k\|_1 + \thalf (\bar\tau_k L + \Gamma) \bar\alpha_k^2 \|\dbar_k\|_2^2 \\
      \leq&\ \bar\alpha_k (\bar\tau_k g_k^T d_k - \|c_k\|_1) + 2 \bar\alpha_k \|c_k\|_1 + \thalf (\bar\tau_k L + \Gamma) \bar\alpha_k^2 \|\dbar_k\|_2^2 + \bar\alpha_k \bar\tau_k g_k^T (\dbar_k - d_k) \\
      \leq&\ -\bar\alpha_k \Delta q(x_k,\bar\tau_k,g_k,H_k,d_k) + 2 \bar\alpha_k \|c_k\|_1 + \thalf (\bar\tau_k L + \Gamma) \bar\alpha_k^2 \|\dbar_k\|_2^2 + \bar\alpha_k \bar\tau_k g_k^T (\dbar_k - d_k).
    \ealigned
  \eequationNN
  Using this inequality, let us now consider two subcases.  (Since the lemma requires $1 \geq \tfrac{\beta_k \bar\xi_k \bar\tau_k}{\bar\tau_k L + \Gamma}$ for all $k \in \N{}$, it is not possible that $\bar\alpha_k = \tfrac{\beta_k \bar\xi_k \bar\tau_k}{\bar\tau_k L + \Gamma}$ in this case.)
  
  \noindent
  \textbf{Case 3a:} If $\bar\alpha_k = \tfrac{\beta_k \Delta q(x_k,\bar\tau_k,\gbar_k,H_k,\dbar_k) - 4 \|c_k\|_1}{(\bar\tau_k L + \Gamma) \|\dbar_k\|_2^2}$, then
  \bequationNN
    \baligned
      &\ \phi(x_k + \bar\alpha_k \dbar_k, \bar\tau_k) - \phi(x_k, \bar\tau_k) \\
      \leq&\ -\bar\alpha_k \Delta q(x_k,\bar\tau_k,g_k,H_k,d_k) + 2 \bar\alpha_k \|c_k\|_1 \\
      &\quad + \thalf \bar\alpha_k (\bar\tau_k L + \Gamma) \(\tfrac{\beta_k \Delta q(x_k,\bar\tau_k,\gbar_k,H_k,\dbar_k) - 4 \|c_k\|_1}{(\bar\tau_k L + \Gamma) \|\dbar_k\|_2^2}\) \|\dbar_k\|_2^2 + \bar\alpha_k \bar\tau_k g_k^T (\dbar_k - d_k) \\
      =&\ -\bar\alpha_k \Delta q(x_k,\bar\tau_k,g_k,H_k,d_k) + \thalf \bar\alpha_k \beta_k \Delta q(x_k,\bar\tau_k,\gbar_k,H_k,\dbar_k) + \bar\alpha_k \bar\tau_k g_k^T (\dbar_k - d_k).
    \ealigned
  \eequationNN
  
  \noindent
  \textbf{Case 3b:} If $\bar\alpha_k = \tfrac{\beta_k \bar\xi_k \bar\tau_k}{\bar\tau_k L + \Gamma} + \theta \beta_k^2 \leq \tfrac{\beta_k \Delta q(x_k,\bar\tau_k,\gbar_k,H_k,\dbar_k) - 4 \|c_k\|_1}{(\bar\tau_k L + \Gamma) \|\dbar_k\|_2^2}$, then
  \bequationNN
    \baligned
      &\ \phi(x_k + \bar\alpha_k \dbar_k, \bar\tau_k) - \phi(x_k, \bar\tau_k) \\
      \leq&\ -\bar\alpha_k \Delta q(x_k,\bar\tau_k,g_k,H_k,d_k) + 2\bar\alpha_k \|c_k\|_1 \\
      &\quad + \thalf \bar\alpha_k (\bar\tau_k L + \Gamma) \(\tfrac{\beta_k \bar\xi_k \bar\tau_k}{\bar\tau_k L + \Gamma} + \theta\beta_k^2\) \|\dbar_k\|_2^2 + \bar\alpha_k \bar\tau_k g_k^T (\dbar_k - d_k) \\
      \leq&\ -\bar\alpha_k \Delta q(x_k,\bar\tau_k,g_k,H_k,d_k) + 2\bar\alpha_k \|c_k\|_1 \\
      &\quad + \thalf \bar\alpha_k \beta_k (\Delta q(x_k,\bar\tau_k,\gbar_k,H_k,\dbar_k) - 4 \|c_k\|_1/\beta_k) + \bar\alpha_k \bar\tau_k g_k^T (\dbar_k - d_k) \\
      \leq&\ -\bar\alpha_k \Delta q(x_k,\bar\tau_k,g_k,H_k,d_k) + \thalf \bar\alpha_k \beta_k \Delta q(x_k,\bar\tau_k,\gbar_k,H_k,\dbar_k) + \bar\alpha_k \bar\tau_k g_k^T (\dbar_k - d_k).
    \ealigned
  \eequationNN
  
  The result follows by combining the conclusions of all cases and subcases.
\eproof

Our next two lemmas provide useful relationships between deterministic (i.e., dependent on $g_k$) and stochastic (i.e., dependent on $\gbar_k$) quantities conditioned on the event that the algorithm has reached $x_k$ as the $k$th iterate.

\blemma\label{lem.expectation}
  For all $k \in \N{}$, $\E_k[\dbar_k] = d_k$, $\E_k[\ubar_k] = u_k$, and $\E_k[\ybar_k] = y_k$.  Moreover, there exists $\kappa_d \in \R{}_{>0}$, independent of $k$ and any run of the algorithm, with
  \bequationNN
    \E_k[\|\dbar_k - d_k\|_2] \leq \kappa_d \sqrt{M}.
  \eequationNN
\elemma
\bproof
  The first statement follows from the fact that, conditioned on the $k$th iterate being~$x_k$, the matrix on the left-hand side of \eqref{eq.system_stochastic} is deterministic and, under Assumption~\ref{ass.deterministic}, it is invertible, along with the fact that expectation is a linear operator.  For the second statement, notice that for any realization of $\gbar_k$, it follows that
  \bequationNN
    \bbmatrix \dbar_k - d_k \\ \ybar_k - y_k \ebmatrix = -\bbmatrix H_k & J_k^T \\ J_k & 0 \ebmatrix^{-1} \bbmatrix \gbar_k - g_k \\ 0 \ebmatrix \implies \|\dbar_k - d_k\|_2 \leq \kappa_d \|\gbar_k - g_k\|_2,
  \eequationNN
  where $\kappa_d \in \R{}_{>0}$ is an upper bound on the norm of the matrix shown above, the existence of which, and independence from $k$, follows under Assumption~\ref{ass.deterministic}.  It also follows from Jensen's inequality, concavity of the square root, and Assumption~\ref{ass.g} that
  \bequationNN
    \E_k[\|\gbar_k - g_k\|_2] \leq \sqrt{\E_k[\|\gbar_k - g_k\|_2^2]} \leq \sqrt{M}.
  \eequationNN
  Combined with the displayed inequality above, the desired conclusion follows.
\eproof

Relationships between inner products involving deterministic and stochastic quantities are the subject of the next lemma.

\blemma\label{lem.product_bounds}
  For all $k \in \N{}$, it follows that
  \bequationNN
    g_k^Td_k \geq \E_k[\gbar_k^T\dbar_k] \geq g_k^Td_k - \zeta^{-1}M\ \ \text{and}\ \ d_k^TH_kd_k \leq \E_k[\dbar_k^TH_k\dbar_k].
  \eequationNN
\elemma
\bproof
  From the first block equation in \eqref{eq.system_stochastic}, it follows that
  \bequationNN
    \baligned
      && H_k(Z_k\wbar_k + v_k) + J_k^T\ybar_k &= -\gbar_k \\
      \iff && Z_k^TH_kZ_k\wbar_k &= -Z_k^T(\gbar_k + H_kv_k) \\
      \iff && Z_k\wbar_k &= -Z_k(Z_k^TH_kZ_k)^{-1}Z_k^T(\gbar_k + H_kv_k),
    \ealigned
  \eequationNN
  from which it follows that
  \bequation\label{eq.bigD1}
    \gbar_k^T\ubar_k = \gbar_k^TZ_k\wbar_k =  -\gbar_k^TZ_k(Z_k^TH_kZ_k)^{-1}Z_k^T(\gbar_k + H_kv_k).
  \eequation
  Following the same line of argument for \eqref{eq.system_deterministic}, it follows that
  \bequation\label{eq.bigD2}
    g_k^Tu_k = -g_k^TZ_k(Z_k^TH_kZ_k)^{-1}Z_k^T(g_k + H_kv_k).
  \eequation
  At the same time, under Assumptions~\ref{ass.H} and \ref{ass.g}, one finds that
  \bequation\label{eq.daniels_special_equation}
    \zeta^{-1}M \geq \E_k[\|Z_k^T(\gbar_k - g_k)\|_{(Z_k^TH_kZ_k)^{-1}}^2] \geq 0.
  \eequation
  One finds that the middle term in this expression can be written as
  \bequationNN
    \baligned
      &\ \E_k[\|Z_k^T(\gbar_k - g_k)\|_{(Z_k^TH_kZ_k)^{-1}}^2] \\
      =&\ \E_k[\|Z_k^T\gbar_k\|_{(Z_k^TH_kZ_k)^{-1}}^2] - 2\E_k[\gbar_k^TZ_k(Z_k^TH_kZ_k)^{-1}Z_k^Tg_k] + \|Z_k^Tg_k\|_{(Z_k^TH_kZ_k)^{-1}}^2 \\
      =&\ \E_k[\|Z_k^T\gbar_k\|_{(Z_k^TH_kZ_k)^{-1}}^2] - \|Z_k^Tg_k\|_{(Z_k^TH_kZ_k)^{-1}}^2.
    \ealigned
  \eequationNN
  Hence, combining \eqref{eq.bigD1}, \eqref{eq.bigD2}, \eqref{eq.daniels_special_equation}, and the fact that $\E_k[\gbar_k] = g_k$ one finds
  \bequationNN
    \baligned
      g_k^Tu_k - \E_k[\gbar_k^T\ubar_k] 
      &= -g_k^TZ_k(Z_k^TH_kZ_k)^{-1}Z_k^T(g_k + H_kv_k) \\
      &\qquad + \E_k[\gbar_k^TZ_k(Z_k^TH_kZ_k)^{-1}Z_k^T(\gbar_k + H_kv_k)] \\
      &= -\|Z_k^Tg_k\|_{(Z_k^TH_kZ_k)^{-1}}^2 + \E_k[\|Z_k^T\gbar_k\|_{(Z_k^TH_kZ_k)^{-1}}^2] \in [0,\zeta^{-1}M].
    \ealigned
  \eequationNN
  The first desired result follows from this fact, $\E_k[\gbar_k^Tv_k] = g_k^Tv_k$, and
  \bequationNN
    g_k^Td_k - \E_k[\gbar_k^T\dbar_k] = g_k^Tu_k + g_k^Tv_k - \E_k[\gbar_k^T\ubar_k + \gbar_k^Tv_k] = g_k^Tu_k - \E_k[\gbar_k^T\ubar_k].
  \eequationNN
  
  Now let us prove the second desired conclusion.  From \eqref{eq.system_stochastic}, it follows that
  \bequationNN
    \baligned
      && H_k(\ubar_k + v_k) &= -\gbar_k - J_k^T\ybar_k \\
      \iff && (\ubar_k + v_k)^TH_k(\ubar_k + v_k) &= -\gbar_k^T(\ubar_k + v_k) - \ybar_k^TJ_k\dbar_k \\
      &&  &= -\gbar_k^T(\ubar_k + v_k) + \ybar_k^Tc_k.
    \ealigned
  \eequationNN
  Following the same argument for \eqref{eq.system_deterministic}, it follows that
  \bequationNN
    (u_k + v_k)^TH_k(u_k + v_k) = -g_k^T(u_k + v_k) + y_k^Tc_k.
  \eequationNN
  Combining these facts, it follows that
  \bequationNN
    \baligned
      &\ \ubar_k^TH_k\ubar_k + 2\ubar_k^TH_kv_k - u_k^TH_ku_k - 2u_k^TH_kv_k \\
      =&\ -\gbar_k^T(\ubar_k + v_k) + g_k^T(u_k + v_k) + (\ybar_k - y_k)^Tc_k,
    \ealigned
  \eequationNN
  which after taking conditional expectation and using Lemma~\ref{lem.expectation} yields
  \bequationNN
    \E_k[\ubar_k^TH_k\ubar_k] - u_k^TH_ku_k = -\E_k[\gbar_k^T\ubar_k] + g_k^Tu_k.
  \eequationNN
  The desired conclusion now follows since
  \bequationNN
    \baligned
      \E_k[\dbar_k^TH_k\dbar_k] - d_k^TH_kd_k
        &= \E_k[(\ubar_k + v_k)^TH_k(\ubar_k + v_k)] - (u_k + v_k)^TH_k(u_k + v_k) \\
        &= \E_k[\ubar_k^TH_k\ubar_k] - u_k^TH_ku_k,
    \ealigned
  \eequationNN
  where again we have used the result of Lemma~\ref{lem.expectation}.
\eproof

In the remainder of our convergence analysis, we consider three cases depending on the behavior of the sequence $\{\bar\tau_k\}$ in a run of the algorithm.  In the deterministic setting, it was proved that the merit parameter sequence eventually remains constant at a value that is sufficiently small to ensure that a primal-dual stationarity measure vanishes (see Lemma~\ref{lem.tau_bound}).  However, under only Assumption~\ref{ass.g}, it is not possible to prove that such behavior is guaranteed for any possible run of Algorithm~\ref{alg.sqp_stochastic}.  Our analysis considers three mutually exclusive and exhaustive events: event $E_{\tau,small}$ that the merit parameter sequence eventually remains constant at a sufficiently small positive value; event $E_{\tau,0}$ that the merit parameter sequence vanishes; and event $E_{\tau,big}$ that the merit parameter sequence eventually remains constant, but at a value that is not sufficiently small.  Under modest assumptions, we prove that $E_{\tau,big}$ occurs with probability zero, and under slightly stronger, but reasonably pragmatic assumptions, we prove that event $E_{\tau,0}$ either does not occur or only occurs in extreme circumstances (e.g., divergence in norm of a subsequence of the stochastic gradient estimates).  This leaves event $E_{\tau,small}$, which we consider first and show that, conditioned on this event, convergence comparable to the deterministic setting is achieved in expectation.

\subsubsection{Constant, Sufficiently Small Merit Parameter}\label{sec.constant_tau_small}

Let us first consider the behavior of the algorithm conditioned on the event that the merit parameter sequence eventually remains constant at a sufficiently small value.  In particular, recalling Lemma~\ref{lem.xi}, let us now make the following assumption.

\bassumption\label{ass.g_conditioned}
  Event $E_{\tau,small}$ occurs in the sense that there exists an iteration number $\kbar_{\tau,\xi} \in \N{}$ and a merit parameter value $\bar\tau_{\min} \in \R{}_{>0}$ such that
  \bequation\label{eq.tau_small}
    \bar\tau_k = \bar\tau_{\min} \leq \tau_k^{trial}\ \ \text{and}\ \ \bar\xi_k = \bar\xi_{\min}\ \ \text{for all}\ \ k \geq \kbar_{\tau,\xi}.
  \eequation
  In addition, the stochastic gradient sequence $\{\gbar_k\}_{k \geq \kbar_{\tau,\xi}}$ satisfies
  \bequationNN
    \E_{k,\tau,small}[\gbar_k] = g_k\ \ \text{and}\ \ \E_{k,\tau,small}[\|\gbar_k - g_k\|_2^2] \leq M,
  \eequationNN
  where $\E_{k,\tau,small}$ denotes expectation with respect to the distribution of $\omega$ conditioned on the event that $E_{\tau,small}$ occurs and the algorithm has reached $x_k$ in iteration $k \in \N{}$.
\eassumption

The inequality $\bar\tau_k \leq \tau_k^{trial}$ in \eqref{eq.tau_small} is critical since it ensures that the model reduction value $\Delta q(x_k,\bar\tau_{\min},g_k,H_k,d_k)$ satisfies the result of Lemma~\ref{lem.Psi} for all $k \geq \kbar_{\tau,\xi}$ with $\bar\tau_{\min}$ in place of $\tau_k$.  In other words, it means that the merit parameter has become small enough such that, if one were to compute the \emph{deterministic} search direction $d_k$ using the \emph{true} gradient $g_k$ at $x_k$, then one would find that it is a direction of sufficient descent for the merit function $\phi(\cdot,\bar\tau_{\min})$ at $x_k$.  The importance of this becomes clear in our final results at the end of this part of our analysis.  The latter part of the assumption reaffirms the properties of the stochastic gradient estimates stated in Assumption~\ref{ass.g}, now conditioned on the occurrence of $E_{\tau,small}$.  With this assumption, the results of Lemma~\ref{lem.expectation} and \ref{lem.product_bounds} continue to hold.  For the sake of brevity, for the rest of this part of our analysis (\S\ref{sec.constant_tau_small}), let us redefine $\E_k[\ \cdot\ ] \equiv \E_{k,\tau,small}[\ \cdot\ ]$.

To derive our main result for this case, our goal is to prove upper bounds in expectation for the positive terms on the right-hand side of the conclusion of Lemma~\ref{lem.key_decrease}.  Let us first consider the last term, which is addressed in our next lemma.

\blemma\label{lem.alphagd}
  Suppose that Assumption~\ref{ass.g_conditioned} holds.  Let $\kappa_g \in \R{}_{>0}$ be an upper bound for $\{\|g_k\|_2\}$, the existence of which follows under Assumption~\ref{ass.deterministic}.  It follows, with $\kappa_d \in \R{}_{>0}$ from Lemma~\ref{lem.expectation} and any $k \geq \kbar_{\tau,\xi}$, that
  \bequationNN
    \E_k[\bar\alpha_k \bar\tau_k g_k^T (\dbar_k - d_k)] \leq \beta_k^2 \theta \bar\tau_{\min} \kappa_g \kappa_d \sqrt{M}.
  \eequationNN
\elemma
\bproof
  For all $k \geq \kbar_{\tau,\xi}$, let $E_k$ be the event that $g_k^T (\dbar_k - d_k) \geq 0$ and let $E_k^c$ be the event that $g_k^T (\dbar_k - d_k) < 0$.  Let $\P_k[\cdot]$ denote probability conditioned on the event that $E_{\tau,small}$ occurs and the algorithm has reached $x_k$ in iteration $k$.  By the Law of Total Expectation, \eqref{eq.tau_small},  and Lemma~\ref{lem.beta_control}, it follows for $k \geq \kbar_{\tau,\xi}$ that
  \bequationNN
    \baligned
      &\ \E_k[\bar\alpha_k \bar\tau_k g_k^T (\dbar_k - d_k)] \\
      =&\ \E_k[\bar\alpha_k \bar\tau_{\min} g_k^T (\dbar_k - d_k) | E_k] \P_k[E_k] + \E_k[\bar\alpha_k \bar\tau_{\min} g_k^T (\dbar_k - d_k) | E_k^c] \P_k[E_k^c] \\
      \leq&\ \bar\alpha_{k,\max} \bar\tau_{\min} \E_k[g_k^T (\dbar_k - d_k) | E_k] \P_k[E_k] + \bar\alpha_{k,\min} \bar\tau_{\min} \E_k[ g_k^T (\dbar_k - d_k) | E_k^c] \P_k[E_k^c].
    \ealigned
  \eequationNN
  By Lemma~\ref{lem.expectation}, this means on one hand that
  \bequationNN
    \baligned
      &\ \E_k[\bar\alpha_k \bar\tau_k g_k^T (\dbar_k - d_k)] \\
      \leq&\ \bar\alpha_{k,\min} \bar\tau_{\min} \E_k[g_k^T (\dbar_k - d_k) | E_k] \P_k[E_k] + \bar\alpha_{k,\min} \bar\tau_{\min} \E_k[ g_k^T (\dbar_k - d_k) | E_k^c] \P_k[E_k^c] \\
      +&\ (\bar\alpha_{k,\max} - \bar\alpha_{k,\min}) \bar\tau_{\min} \E_k[g_k^T (\dbar_k - d_k) | E_k] \P_k[E_k] \\
      =&\ (\bar\alpha_{k,\max} - \bar\alpha_{k,\min}) \bar\tau_{\min} \E_k[g_k^T (\dbar_k - d_k) | E_k] \P_k[E_k],
    \ealigned
  \eequationNN
  while on the other hand that
  \bequationNN
    \baligned
      &\ \E_k[\bar\alpha_k \bar\tau_k g_k^T (\dbar_k - d_k)] \\
      \leq&\ \bar\alpha_{k,\max} \bar\tau_{\min} \E_k[g_k^T (\dbar_k - d_k) | E_k] \P_k[E_k] + \bar\alpha_{k,\max} \bar\tau_{\min} \E_k[ g_k^T (\dbar_k - d_k) | E_k^c] \P_k[E_k^c] \\
      +&\ (\bar\alpha_{k,\min} - \bar\alpha_{k,\max}) \bar\tau_{\min} \E_k[g_k^T (\dbar_k - d_k) | E_k^c] \P_k[E_k^c] \\
      =&\ (\bar\alpha_{k,\min} - \bar\alpha_{k,\max}) \bar\tau_{\min} \E_k[g_k^T (\dbar_k - d_k) | E_k^c] \P_k[E_k^c]
    \ealigned
  \eequationNN
  Combining these facts, it follows that
  \bequationNN
    \baligned
      &\ \E_k[\bar\alpha_k \bar\tau_k g_k^T (\dbar_k - d_k)] \\
      \leq&\ \thalf (\bar\alpha_{k,\max} - \bar\alpha_{k,\min}) \bar\tau_{\min} (\E_k[g_k^T (\dbar_k - d_k) | E_k] \P_k[E_k] - \E_k[g_k^T (\dbar_k - d_k) | E_k^c] \P_k[E_k^c]).
    \ealigned
  \eequationNN
  Observe that, by the Law of Total Expectation, it follows that
  \bequationNN
    \baligned
      \E_k[g_k^T (\dbar_k - d_k) | E_k] \P_k[E_k] 
      \leq&\ \E_k[\|g_k\|_2 \|\dbar_k - d_k\|_2 | E_k] \P_k[E_k] \\
      =&\ \E_k[\|g_k\|_2 \|\dbar_k - d_k\|_2] - \E_k[\|g_k\|_2 \|\dbar_k - d_k\|_2 | E_k^c] \P_k[E_k^c] \\
      \leq&\ \|g_k\|_2 \E_k[\|\dbar_k - d_k\|_2],
    \ealigned
  \eequationNN
  and, in a similar manner,
  \bequationNN
    \baligned
      -\E_k[g_k^T (\dbar_k - d_k) | E_k^c] \P_k[E_k^c]
      \leq&\ \E_k[\|g_k\|_2 \|\dbar_k - d_k\|_2 | E_k^c] \P_k[E_k^c] \\
      =&\ \E_k[\|g_k\|_2 \|\dbar_k - d_k\|_2] - \E_k[\|g_k\|_2 \|\dbar_k - d_k\|_2 | E_k] \P_k[E_k] \\
      \leq&\ \|g_k\|_2 \E_k[\|\dbar_k - d_k\|_2].
    \ealigned
  \eequationNN
  Combining these results with Lemma~\ref{lem.beta_control} and Lemma~\ref{lem.expectation} yield the result.
\eproof

Our next result addresses the middle term on the right-hand side of Lemma~\ref{lem.key_decrease}.

\blemma\label{lem.Deltaq}
  Suppose Assumption~\ref{ass.g_conditioned} holds.  Then, for all $k \geq \kbar_{\tau,\xi}$, it follows that
  \bequationNN
    \E_k[\Delta q(x_k,\bar\tau_k,\gbar_k,H_k,\dbar_k)] \leq \Delta q(x_k,\bar\tau_{\min},g_k,H_k,d_k) + \bar\tau_{\min} \zeta^{-1} M.
  \eequationNN
\elemma
\bproof
  Consider arbitrary $k \geq \kbar_{\tau,\xi}$.  From \eqref{def.merit_model_reduction}, \eqref{eq.tau_small},  Lemma~\ref{lem.product_bounds}, Jensen's inequality, and convexity of $\max\{\cdot,0\}$, it follows that
  \bequationNN
    \baligned
      \E_k[\Delta q(x_k,\bar\tau_k,\gbar_k,H_k,\dbar_k)]
      &= \E_k[-\bar\tau_{\min}(\gbar_k^T\dbar_k + \thalf \max \{\dbar_k^TH_k\dbar_k,0\}) + \|c_k\|_1] \\
      &\leq -\bar\tau_{\min}(g_k^Td_k + \thalf \max\{d_k^TH_kd_k,0\}) + \bar\tau_{\min} \zeta^{-1}M + \|c_k\|_1 \\
      &= \Delta q(x_k,\bar\tau_{\min},g_k,H_k,d_k) + \bar\tau_{\min} \zeta^{-1} M,
    \ealigned
  \eequationNN
  as desired.
\eproof

We now prove our main theorem for this part of our analysis, where we define
\bequationNN
  \E_{\tau,small} [\ \cdot\ ] = \E[\ \cdot\ |\ \text{Assumption~\ref{ass.g_conditioned} holds}\ ].
\eequationNN
The theorem considers the behavior of a certain sequence of model reduction values, and the subsequent corollary translates the result of the theorem to the behavior of the sequences of constraint violations and stationarity measures.

\btheorem\label{th.stochastic_tau_constant_small}
  Suppose that Assumption~\ref{ass.g_conditioned} holds and the sequence $\{\beta_k\}$ is chosen such that $\beta_k \bar\xi_k \bar\tau_k/(\bar\tau_k L + \Gamma) \in (0,1]$ for all $k \geq \kbar_{\tau,\xi}$.  Define
  \bequationNN
    \Abar := \tfrac{\bar\xi_{\min} \bar\tau_{\min}}{\bar\tau_{\min} L + \Gamma}\ \ \text{and}\ \ \Mbar := \bar\tau_{\min} \big(\thalf (\Abar + \theta)  \zeta^{-1} M + \theta \kappa_g \kappa_d \sqrt{M}\big).
  \eequationNN
  If $\beta_k = \beta \in (0,2\Abar/(\Abar + \theta))$ for all $k \geq \kbar_{\tau,\xi}$, then
  \bequation\label{eq.beta_fixed}
    \baligned
      &\ \E_{\tau,small}\left[\tfrac{1}{k+1} \sum_{j=\kbar_{\tau,\xi}}^{\kbar_{\tau,\xi}+k} \Delta q(x_j,\bar\tau_{\min},g_j,H_j,d_j) \right] \\
      \leq&\ \tfrac{\beta \Mbar}{\Abar - \thalf(\Abar + \theta)\beta} + \tfrac{\E_{\tau,small}[\phi(x_{\kbar_{\tau,\xi}},\bar\tau_{\min})] - \phi_{\min}}{(k+1) \beta(\Abar - \thalf(\Abar + \theta)\beta)} \xrightarrow{k\to\infty} \tfrac{\beta \Mbar}{\Abar - \thalf(\Abar + \theta)\beta},
    \ealigned
  \eequation
  where $\phi_{\min} \in \R{}_{>0}$ is a lower bound for $\phi(\cdot,\bar\tau_{\min})$ over $\Xcal$, the existence of which follows by Assumption~\ref{ass.deterministic}.  On the other hand, if $\sum_{k=\kbar_{\tau,\xi}}^\infty \beta_k = \infty$ and $\sum_{k=\kbar_{\tau,\xi}}^\infty \beta_k^2 < \infty$, then
  \bequation\label{eq.beta_diminishing}
    \lim_{k \to \infty} \E_{\tau,small}\left[ \tfrac{1}{\(\sum_{j=\kbar_{\tau,\xi}}^{\kbar_{\tau,\xi}+k} \beta_j\)} \sum_{j=\kbar_{\tau,\xi}}^{\kbar_{\tau,\xi}+k} \beta_j\Delta q(x_j,\bar\tau_{\min},g_j,H_j,d_j) \right] = 0.
  \eequation
\etheorem
\bproof
  Consider arbitrary $k \geq \kbar_{\tau,\xi}$.  It follows from the definition of $\Abar$, Lemma~\ref{lem.beta_control}, and the fact that $\beta_k \in (0,1]$ that $\Abar \beta_k \leq \bar\alpha_k \leq (\Abar + \theta) \beta_k$.  Hence, it follows from Lemmas~\ref{lem.deterministic_to_stochastic}(d), \ref{lem.key_decrease}, \ref{lem.alphagd}, and \ref{lem.Deltaq} that, under the conditions of the theorem,
  \bequationNN
    \baligned
      &\ \E_k[\phi(x_k + \bar\alpha_k \dbar_k, \bar\tau_k)] - \E_k[\phi(x_k, \bar\tau_k)] \\
      \leq&\ \E_k[-\bar\alpha_k \Delta q(x_k,\bar\tau_k,g_k,H_k,d_k) + \thalf \bar\alpha_k \beta_k \Delta q(x_k,\bar\tau_k,\gbar_k,H_k,\dbar_k) + \bar\alpha_k \bar\tau_k g_k^T (\dbar_k - d_k)] \\
      \leq&\ -\beta_k \big(\Abar - \thalf (\Abar + \theta) \beta_k\big) \Delta q(x_k,\bar\tau_{\min},g_k,H_k,d_k) + \beta_k^2 \Mbar.
    \ealigned
  \eequationNN
  For the scenario of $\{\beta_k\}$ being a constant sequence for $k \geq \kbar_{\tau,\xi}$, one finds from above, taking total expectation conditioned on \eqref{eq.tau_small}, that, for all $k \geq \kbar_{\tau,\xi}$,
  \begin{multline*}
    \E_{\tau,small}[\phi(x_k + \bar\alpha_k \dbar_k, \bar\tau_{\min})] - \E_{\tau,small}[\phi(x_k, \bar\tau_{\min})] \\
    \leq - \beta(\Abar - \thalf(\Abar + \theta)\beta) \E_{\tau,small}[\Delta q(x_k,\bar\tau_{\min},g_k,H_k,d_k)] + \beta^2 \Mbar.
  \end{multline*}
  Summing this inequality for $j \in \{\kbar_{\tau,\xi},\dots,\kbar_{\tau,\xi} + k\}$, one finds by Assumption~\ref{ass.deterministic} that
  \bequationNN
    \baligned
      &\ \phi_{\min} - \E_{\tau,small}[\phi(x_{\kbar_{\tau,\xi}},\bar\tau_{\min})] \\
      \leq&\ \E_{\tau,small}[\phi(x_{\kbar_{\tau,\xi}+k+1},\bar\tau_{\min})] - \E_{\tau,small}[\phi(x_{\kbar_{\tau,\xi}},\bar\tau_{\min})] \\
      \leq&\ -\beta(\Abar - \thalf(\Abar + \theta)\beta) \E_{\tau,small}\left[\sum_{j=\kbar_{\tau,\xi}}^{\kbar_{\tau,\xi}+k} \Delta q(x_j,\bar\tau_{\min},g_j,H_j,d_j) \right] + (k+1) \beta^2 \Mbar,
    \ealigned
  \eequationNN
  from which \eqref{eq.beta_fixed} follows.  Now consider the scenario of $\{\beta_k\}$ diminishing as described.  It follows that for sufficiently large $k \geq \kbar_{\tau,\xi}$ one finds $\beta_k \leq \Abar/(\Abar + \theta)$; hence, let us assume without loss of generality that, for all $k \geq \kbar_{\tau,\xi}$, one has $\beta_k \leq \Abar/(\Abar + \theta)$, which implies $\Abar - \thalf(\Abar + \theta)\beta_k \geq \thalf \Abar$.  Similar to above, it follows for all $k \geq \kbar_{\tau,\xi}$ that
  \begin{multline*}
    \E_{\tau,small}[\phi(x_k + \bar\alpha_k \dbar_k, \bar\tau_{\min})] - \E_{\tau,small}[\phi(x_k, \bar\tau_{\min})] \\ \leq - \thalf \Abar \beta_k \E_{\tau,small}[\Delta q(x_k,\bar\tau_{\min},g_k,H_k,d_k)] + \beta_k^2 \Mbar.
  \end{multline*}
  Summing this inequality for $j \in \{\kbar_{\tau,\xi},\dots,\kbar_{\tau,\xi} + k\}$, one finds by Assumption~\ref{ass.deterministic} that
  \bequationNN
    \baligned
      &\ \phi_{\min} - \E_{\tau,small}[\phi(x_{\kbar_{\tau,\xi}},\bar\tau_{\min})] \\
      \leq&\ \E_{\tau,small}[\phi(x_{\kbar_{\tau,\xi}+k+1},\bar\tau_{\min})] - \E_{\tau,small}[\phi(x_{\kbar_{\tau,\xi}},\bar\tau_{\min})] \\
      \leq&\ -\thalf \Abar \E_{\tau,small}\left[ \sum_{j=\kbar_{\tau,\xi}}^{\kbar_{\tau,\xi}+k} \beta_j \Delta q(x_j,\bar\tau_{\min},g_j,H_j,d_j)\right] + \Mbar \sum_{j=\kbar_{\tau,\xi}}^{\kbar_{\tau,\xi}+k} \beta_j^2.
    \ealigned
  \eequationNN
  Rearranging this inequality yields
  \begin{multline*}
    \E_{\tau,small} \left[ \sum_{j=\kbar_{\tau,\xi}}^{\kbar_{\tau,\xi}+k} \beta_j \Delta q(x_j,\bar\tau_{\min},g_j,H_j,d_j)\right] \\ \leq \tfrac{2(\E_{\tau,small}[\phi(x_{\kbar_{\tau,\xi}},\bar\tau_{\min})] - \phi_{\min})}{\Abar} + \tfrac{2\Mbar}{\Abar} \sum_{j=\kbar_{\tau,\xi}}^{\kbar_{\tau,\xi}+k} \beta_j^2,
  \end{multline*}
  from which \eqref{eq.beta_diminishing} follows.
\eproof

\bcorollary\label{cor.stochastic_tau_finite_large} 
  Under the conditions of Theorem~\ref{th.stochastic_tau_constant_small}, the following hold true.
  \benumerate
    \item[(a)] If $\beta_k = \beta \in (0,2\Abar/(\Abar + \theta))$ for all $k \geq \kbar_{\tau,\xi}$, then
    \bequationNN
      \baligned
        &\ \E_{\tau,small}\left[\tfrac{1}{k+1} \sum_{j=\kbar_{\tau,\xi}}^{\kbar_{\tau,\xi}+k} \(\tfrac{\|g_j + J_j^Ty_j\|_2^2}{\kappa_H^2} + \|c_j\|_2\) \right] \xrightarrow{k\to\infty} \tfrac{2\kappa_\Psi \beta \Mbar}{\kappa_q \bar\tau_{\min} (\Abar - \thalf(\Abar + \theta)\beta)}.
      \ealigned
    \eequationNN
    \item[(b)] If $\sum_{k=\kbar_{\tau,\xi}}^\infty \beta_k = \infty$ and $\sum_{k=\kbar_{\tau,\xi}}^\infty \beta_k^2 < \infty$, then
    \bequationNN
       \lim_{k\to\infty}\E_{\tau,small}\left[ \tfrac{1}{\(\sum_{j=\kbar_{\tau,\xi}}^{\kbar_{\tau,\xi}+k} \beta_j\)} \sum_{j=\kbar_{\tau,\xi}}^{\kbar_{\tau,\xi}+k} \beta_j\(\tfrac{\|g_j + J_j^Ty_j\|_2^2}{\kappa_H^2} + \|c_j\|_2\) \right] = 0,
    \eequationNN
    from which it follows that \bequationNN
      \liminf_{k\to\infty}\  \E_{\tau,small} [\kappa_H^{-2} \|g_k + J_k^Ty_k\|_2^2 + \|c_k\|_2] = 0.
    \eequationNN
  \eenumerate
  In addition, in either case, there exists $\delta_x \in \R{}_{>0}$ such that if $\|x_k - x_*\|_2 \leq \delta_x$ for some stationary point $(x_*,y_*) \in \R{n} \times \R{m}$ for \eqref{prob.f_nonlinear_stochastic}, then for any $\delta_g \in \R{}_{>0}$ one finds
  \bequationNN
    \left\|\bbmatrix\gbar_k - \nabla f(x_*) \\ c_k \ebmatrix\right\|_2 \leq \delta_g\ \ \implies\ \ \|\ybar_k - y_*\|_2 \leq 2\delta_g.
  \eequationNN
\ecorollary
\bproof
  Parts (a) and (b) follow by combining the results of Lemmas~\ref{lem.Psi_1} and \ref{lem.Psi}, the relation \eqref{eq.gJy}, and Theorem~\ref{th.stochastic_tau_constant_small}.  The remainder follows with Lemma~\ref{lem.stationary} since for $x_k$ sufficiently close to $x_*$, one obtains with $g_* := \nabla f(x_*)$ and $c_* := c(x_*) = 0$ that
  \bequationNN
    \|\ybar_k - y_*\|_2 \leq \left\|\bbmatrix H_k & J_k^T \\ J_k & 0 \ebmatrix^{-1} \bbmatrix \gbar_k \\ c_k \ebmatrix - \bbmatrix H_k & J_*^T \\ J_* & 0 \ebmatrix^{-1} \bbmatrix g_* \\ c_* \ebmatrix \right\|_2 \leq 2 \left\|\bbmatrix\gbar_k - g_* \\ c_k\ebmatrix\right\|_2,
  \eequationNN
  from which the desired conclusion follows.
\eproof

We close our analysis of this case with the following remark.

\begin{remark}
  Consideration of the conclusion of Corollary~\ref{cor.stochastic_tau_finite_large}(a) reveals the close relationship between our result and a conclusion that one reaches for a stochastic (sub)gradient method in an unconstrained setting.  Notice that
  \bequationNN
    \tfrac{2\kappa_\Psi \beta \Mbar}{\kappa_q \bar\tau_{\min} (\Abar - \thalf(\Abar + \theta)\beta)} = \tfrac{2\kappa_\Psi \beta \bar\tau_{\min} \(\thalf \(\tfrac{\bar\xi_{\min} \bar\tau_{\min}}{\bar\tau_{\min} L + \Gamma} + \theta\) \zeta^{-1} M + \theta \kappa_g \kappa_d \sqrt{M}\)}{\kappa_q \bar\tau_{\min} \(\tfrac{\bar\xi_{\min} \bar\tau_{\min}}{\bar\tau_{\min} L + \Gamma} - \thalf\(\tfrac{\bar\xi_{\min} \bar\tau_{\min}}{\bar\tau_{\min} L + \Gamma} + \theta\)\beta\)}.
  \eequationNN
  Our first observation is a common one for the unconstrained setting: The value above is directly proportional to~$\beta$.  To reduce this value, one should choose smaller~$\beta$, but the downside of choosing smaller~$\beta$ is that the algorithm takes shorter steps, meaning that it takes longer for this limiting value to be approached $($recall~\eqref{eq.beta_fixed}$)$.  On the other hand, while larger $\beta$ means that the algorithm takes larger steps, this comes at the cost of a larger limiting value.  A second observation, unique for our algorithm, is the influence of $\theta$.  The quantity above is directly proportional to $\theta$, meaning that the optimal choice in terms of reducing this value is $\theta=0$, in which case one obtains
  \bequationNN
    \tfrac{2\kappa_\Psi \beta \Mbar}{\kappa_q \bar\tau_{\min} (\Abar - \thalf(\Abar + \theta)\beta)} \xrightarrow{\theta \to 0} \tfrac{\kappa_\Psi \beta \zeta^{-1} M}{\kappa_q(1 - \thalf \beta)}.
  \eequationNN
  However, this results in a non-adaptive algorithm with $\bar\alpha_k = \beta_k\bar\xi_k\bar\tau_k/(\bar\tau_k L + \Gamma)$ for all $k \in \N{}$.  This choice has some theoretical benefits $($see also our discussion in \S\ref{sec.conclusion}$)$, but we have found this conservative choice to be detrimental in practice.
\end{remark}

\subsubsection{Poor Merit Parameter Behavior}

Theorem~\ref{th.stochastic_tau_constant_small} and Corollary~\ref{cor.stochastic_tau_finite_large} show desirable convergence properties in expectation of Algorithm~\ref{alg.sqp_stochastic} in the event that the merit parameter sequence eventually remains constant at a value that is sufficiently small.  This captures behavior similar to that of Algorithm~\ref{alg.sqp_adaptive} in the deterministic setting, in which the merit parameter is \emph{guaranteed} to behave in this manner.  However, for the stochastic  Algorithm~\ref{alg.sqp_stochastic}, one of two other events are possible, which we now define mathematically as follows:
\bitemize
  \item Event $E_{\tau,big}$: there exists infinite $\overline\Kcal_\tau \subseteq \N{}$ and $\bar\tau_{big} \in \R{}_{>0}$ such that
  \bequationNN
    \bar\tau_k = \bar\tau_{big} > \tau_k^{trial}\ \ \text{and}\ \ \bar\xi_k = \bar\xi_{\min}\ \ \text{for all}\ \ k \in \overline\Kcal_\tau.
  \eequationNN
  Since $\bar\tau_k^{trial} \geq \bar\tau_k$ for all $k \in \N{}$, this means $\bar\tau_k^{trial} > \tau_k^{trial}$ for all $k \in \overline\Kcal_\tau$.
  \item Event $E_{\tau,0}$: $\{\bar\tau_k\} \searrow 0$.
\eitemize
Our goal in this part of our analysis is to argue that these events, exhibiting what we refer to as poor behavior of the merit parameter sequence, are either impossible or can only occur in extreme circumstances in practice.  For these considerations, let us return to assume that Assumption~\ref{ass.g} (not Assumption~\ref{ass.g_conditioned}) holds.

Let us first consider event $E_{\tau,big}$.  We show under a modest assumption that this event occurs with probability zero, which is to say that the merit parameter eventually becomes sufficiently small with probability one.  As shown above in the definition of $E_{\tau,big}$, the merit parameter remaining too large requires that the stochastic trial value $\bar\tau_k^{trial}$ \emph{consistently} overestimates the deterministic trial value $\tau_k^{trial}$.  The following proposition shows that under a modest assumption about the behavior of the stochastic gradients and corresponding search directions, this behavior occurs with probability zero.  The subsequent proposition then provides an example showing that our modest assumption holds for an archetypal distribution of the stochastic gradients.

\begin{proposition}\label{prop.p}
  If there exists $p \in (0,1]$ such that, for all $k \in \N{}$,
  \bequationNN
    \P_k[\gbar_k^T\dbar_k + \max\{\dbar_k^TH_k\dbar_k,0\} \geq g_k^Td_k + \max\{d_k^TH_kd_k,0\}] \geq p,
  \eequationNN
  then $E_{\tau,big}$ occurs with probability zero.
\end{proposition}
\bproof
  If, in any run of the algorithm, $g_k^Td_k + \max\{d_k^TH_kd_k,0\} \leq 0$ for all sufficiently large $k \in \N{}$, then $\tau_k^{trial} = \infty$ for all sufficiently large $k \in \N{}$ and event $E_{\tau,big}$ does not occur.  Hence, let us define $\Kcal_{gd} \subseteq \N{}$ as the set of indices such that $k \in \Kcal_{gd}$ if and only if $g_k^Td_k + \max\{d_k^TH_kd_k,0\} > 0$, and let us restrict attention to runs in which $\Kcal_{gd}$ is infinite.  For any $k \in \Kcal_{gd}$, it follows that the inequality $\gbar_k^T\dbar_k + \max\{\dbar_k^TH_k\dbar_k,0\} \geq g_k^Td_k + \max\{d_k^TH_kd_k,0\}$ holds if and only if
  \bequationNN
    \bar\tau_k^{trial} = \tfrac{(1-\sigma)\|c_k\|_1}{\gbar_k^T\dbar_k + \max\{\dbar_k^TH_k\dbar_k,0\}} \leq \tfrac{(1-\sigma)\|c_k\|_1}{g_k^Td_k + \max\{d_k^TH_kd_k,0\}} = \tau_k^{trial}.
  \eequationNN
  Hence, it follows from the conditions of the proposition, the fact that $\bar\tau_k \leq \bar\tau_k^{trial}$ for all $k \in \N{}$, and the fact that $\Kcal_{gd}$ is infinite, that for any $k \in \N{}$ the probability is one that for a subsequent iteration number $\khat \geq k$ one finds $\bar\tau_{\khat} \leq \bar\tau_{\khat}^{trial} \leq \tau_{\khat}^{trial}$.  This, the fact that Lemma~\ref{lem.tau_bound} implies that $\{\tau_k^{trial}\}$ is bounded away from zero, and the fact that if the merit parameter is ever decreased then it is done so by a constant factor, shows that one has $\bar\tau_k \leq \tau_k^{trial}$ for all sufficiently large $k \in \N{}$ with probability one.
\eproof

As a concrete example of a setting that offers the minimum probability required in Proposition~\ref{prop.p}, we offer the following.  This is clearly only one of many example situations that one could consider to mimic real-world scenarios.

\begin{example}
  If, for all $k \in \N{}$, one has $H_k \succ 0$ and $\gbar_k \sim \Ncal(g_k,\Sigma_k)$ for some $\Sigma_k \in \mathbb{S}^n$ with $\Sigma_k \succ 0$, then the condition in Proposition~\ref{prop.p} holds with $p = \thalf$.
\end{example}
\bproof
  Let $k \in \N{}$ be arbitrary.  The tangential component of the search direction is $\ubar_k = Z_k\wbar_k$, where, under Assumption~\ref{ass.H} and the stated conditions, $\wbar_k = -(Z_k^TH_kZ_k)^{-1}Z_k^T(\gbar_k + H_kv_k)$.  Plugging in this solution and simplifying yields
  \bequationNN
    \gbar_k^T\dbar_k + \dbar_k^TH_k\dbar_k = v_k^TH_k^{1/2}(I - H_k^{1/2}Z_k(Z_k^TH_kZ_k)^{-1}Z_k^TH_k^{1/2})(H_k^{-1/2}\gbar_k + H_k^{1/2}v_k).
  \eequationNN
  Since $\gbar_k$ is normally distributed with mean $g_k$, it follows that this value is normally distributed with mean of the same form, but with $g_k$ in place of $\gbar_k$ (see, e.g., \cite{Tong12}).  Since a normally distributed random variable takes values greater than or equal to its expected value with probability $\thalf$, the conclusion follows.
\eproof

Let us now consider the event $E_{\tau,0}$.  One can learn from Lemmas~\ref{lem.bound_u} and \ref{lem.tau_bound} from the deterministic setting that the following holds true.

\begin{proposition}\label{lem.tau_bound_stochastic}
  Consider an arbitrary constant $g_{\max} \in \R{}_{>0}$.  If, for a run of Algorithm~\ref{alg.sqp_stochastic}, the stochastic gradient estimates satisfy $\|\gbar_k - g_k\|_2 \leq g_{\max}$ for all $k \in \N{}$, then the sequence of tangential step components $\{\ubar_k\}$ is bounded, and there exists $\kbar_\tau \in \N{}$ and $\bar\tau_{\min} \in \R{}_{>0}$ such that $\bar\tau_k = \bar\tau_{\min}$ for all $k \geq \kbar_\tau$.
\end{proposition}
\bproof
  Boundedness in norm of the tangential step components follows in the same manner as in Lemma~\ref{lem.bound_u} with $(\gbar_k,\ubar_k)$ in place of $(g_k,u_k)$.  Further, the claimed behavior of the merit parameter sequence follows in the same manner as in the proof of Lemma~\ref{lem.tau_bound} using $(\gbar_k,\dbar_k,\ubar_k)$ in place of $(g_k,d_k,u_k)$, where in place of the constants $(\kappa_{\tau,1},\kappa_{\tau,2})$ one derives constants $(\bar\kappa_{\tau,1},\bar\kappa_{\tau,2})$ whose value depends on $g_{\max}$ as well as the upper bound on the sequence $\{\|g_k\|_2\}$ (under Assumption~\ref{ass.deterministic}).
\eproof

By Proposition~\ref{lem.tau_bound_stochastic}, if the differences between the stochastic gradient estimates and true gradients are bounded in norm, then the merit parameter sequence will not vanish, i.e., event $E_{\tau,0}$ will not occur.  This is guaranteed if the distributions defining the stochastic gradients $\{\gbar_k\}$ ensure uniform boundedness or, e.g., if
\bequationNN
  f(x) = \tfrac{1}{N} \sum_{i=1}^N f_i(x)\ \ \text{and}\ \ \gbar_k := \nabla f_{i_k}(x_k)\ \ \text{for all}\ \ k \in \N{},
\eequationNN
where the component functions $\{f_i\}$ have bounded derivatives over a set containing the iterates and in each iteration $i_k$ is randomly sampled uniformly from $\{1,\dots,N\}$.

\section{Numerical Results}\label{sec.numerical}

In this section, we demonstrate the empirical performance of our proposed Algorithm~\ref{alg.sqp_adaptive} (for the deterministic setting) and Algorithm~\ref{alg.sqp_stochastic} (for the stochastic setting) using Matlab implementations.  We consider their performance on a subset of the equality constrained problems from the CUTE collection \cite{BongConnGoulToin95}.  Specifically, of the 123 such problems in the set, we selected those for which $(i)$ $f$ is \emph{not} a constant function, $(ii)$~$n+m \leq 1000$, and $(iii)$ the LICQ held at all iterates in all runs of all algorithms that we ran.  This selection resulted in a total of 49 problems.  Each problem comes with an initial point, which we used in our experiments.

\subsection{Deterministic Setting}

Our goal in this setting is to demonstrate that, in practice, our proposed Algorithm~\ref{alg.sqp_adaptive} (``SQP Adaptive'') is as reliable a method as the state-of-the-art Algorithm~\ref{alg.sqp_line_search} (``SQP Backtracking'').  We do not claim that ``SQP Adaptive'' is as efficient as ``SQP Backtracking'' since, as has been verified by others in the literature, the line search scheme is very effective across a broad range of problems.  That said, since our algorithm for the stochastic setting is based on ``SQP Adaptive,'' it is at least of interest to demonstrate that this approach is as reliable as ``SQP Backtracking'' in practice.  For these experiments, we chose each $H_k$ to be the Hessian of the Lagrangian at $(x_k,y_{k-1})$.  For both algorithms, for any $k$ such that the inertia of the matrix in \eqref{eq.system_deterministic} is not correct with this choice, a multiple of the identity is added in an iterative manner until the correct inertia is attained.  This is a common strategy in state-of-the-art constrained optimization software; see, e.g., \cite{WaecBieg06}.

\begin{figure}[ht]
  \centering
  \includegraphics[width=0.425\textwidth,clip=true,trim=60 10 90 50]{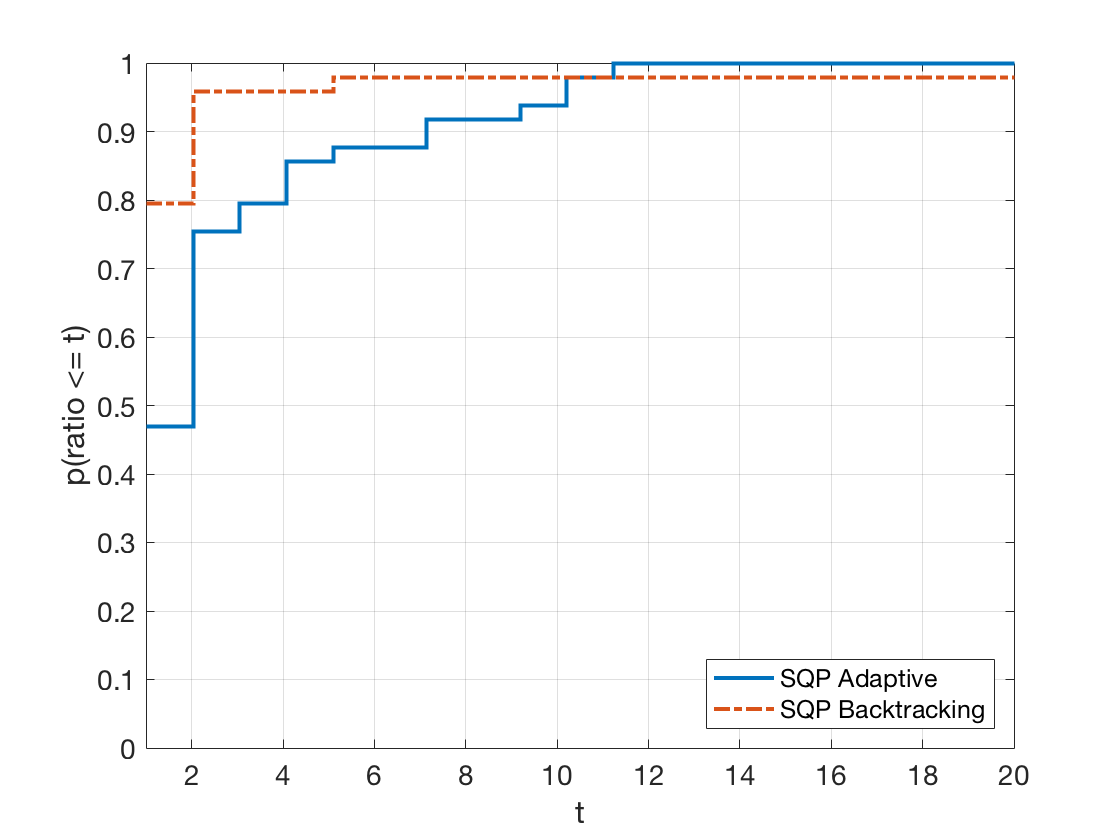} \quad
  \includegraphics[width=0.425\textwidth,clip=true,trim=60 10 90 50]{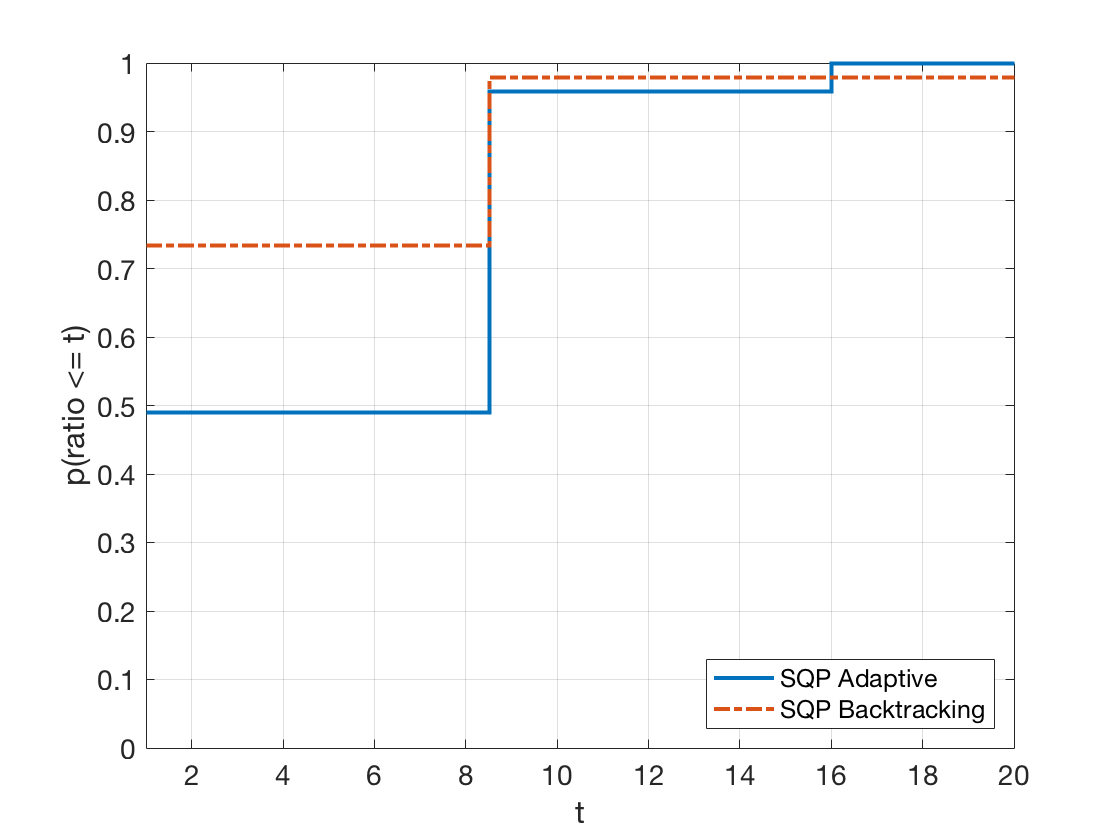}
  \caption{Performance profiles for ``SQP Adaptive'' and ``SQP Backtracking'' for problems from the CUTE test set in terms of iterations (left) and function evaluations (right).}
  \label{fig.perf_deterministic}
\end{figure}

For our experiments, the parameters were set as: $\tau_{-1} = 1$, $\epsilon = 10^{-6}$, $\sigma = 1/2$, $\eta = 10^{-4}$, $\rho = 3$, $L_{-1} = 1$, $\gamma_{-1,i} = 1$, $\nu=1/2$, and $\alpha = 1$.  In Line~\ref{step.initialize} of Algorithm \ref{alg.sqp_adaptive}, all Lipschitz constant estimates were set as $1/2$ times the estimates from the previous iteration.  A run terminated with a message of success if iteration $k \leq 10^4$ yielded
\bequationNN
  \|g_k + J_k^Ty_k\|_\infty \leq 10^{-6} \max\{1,\|g_0 + J_0^Ty_0\|_\infty\}\ \ \text{and}\ \ \|c_k\|_\infty \leq 10^{-6} \max\{1,\|c_0\|_\infty\};
\eequationNN
otherwise, the run was considered a failure.  Figure \ref{fig.perf_deterministic} provides Dolan-Mor\'e performance profiles \cite{DolaMore02} for iterations and function evaluations required by the two methods.  (The profiles are capped at $t = 20$.)  As expected, the performance of ``SQP Backtracking'' was typically better than that of ``SQP Adaptive.''  That said, ``SQP Adaptive'' was as reliable as this state-of-the-art approach.  Over all iterations of all runs of ``SQP Adaptive,'' the stepsize $\alpha_k$ was chosen less than one $40.9\%$ of the time, equal to one $41.8\%$ of the time, and greater than one $17.3\%$ percent of the time.

\subsection{Stochastic Setting}

Our goal in this setting is to compare the performance of our proposed Algorithm \ref{alg.sqp_stochastic} (``Stochastic SQP'') against that of a stochastic subgradient method (``Stochastic subgradient'') applied to minimize the exact penalty function~\eqref{eq.penalty_function} (which represents the current state-of-the-art for constrained stochastic optimization).  For these experiments, we used our test set of 49 CUTE problems, but considered multiple runs for different levels of noise.  In particular, for a given run of an algorithm, we fixed $\epsilon_N \in \{10^{-8}, 10^{-4},10^{-2},10^{-1}\}$, then for each iteration set the stochastic gradient estimate as $\gbar_k = \Ncal(g_k,\epsilon_N I)$.  For each problem and noise level, we ran 10 instances.  This led to a total of 490 problem instances for each algorithm and noise level.  Each run of ``Stochastic SQP'' was given a budget of $1000$ iterations while each run of ``Stochastic Subgradient'' was given a budget of $10000$ iterations.  We tuned the value of $\tau$ individually for each problem instance for ``Stochastic Subgradient.''  In particular, for each problem instance, we ran the algorithm for the $11$ values $\tau\in \{10^{-10},10^{-9},\dots,10^{-1},10^0\}$ and selected the value for that instance that led to the best results in terms of feasibility and optimality errors (see below).  Overall, this means that for each problem, ``Stochastic Subgradient'' was given 110 times the number of iterations that were allowed for ``Stochastic SQP.''  (This broad range of $\tau$ was needed by ``Stochastic Subgradient'' to obtain its best results.  The selected $\tau$ values were roughly evenly distributed over the set from $10^{-10}$ to $10^0$.)

For both methods, the Lipschitz constants $L$ and $\Gamma = \sum_{i=1}^m\gamma_i$ were estimated using differences of gradients near the initial point and kept fixed for all subsequent iterations.  (This process was done so that $L$ and $\Gamma$ were the same for both methods for each problem.)  For ``Stochastic SQP,'' we set $H_k=I$ for all $k$ for fairness of comparison with the (first-order) subgradient method.  The other inputs for ``Stochastic SQP'' were set as: $\bar{\tau}_{-1}=1$, $\epsilon = 10^{-6}$, $\sigma = 1/2$, $\bar{\xi}_{-1}=1$, $\theta = 10$, and $\beta_k = 1$ for all $k$.  ``Stochastic Subgradient'' was run with a constant stepsize $\tfrac{\tau}{\tau L + \Gamma}$ for all $k$.

\begin{figure}[ht]
  \centering
  \includegraphics[width=0.45\textwidth,clip=true,trim=20 5 90 50]{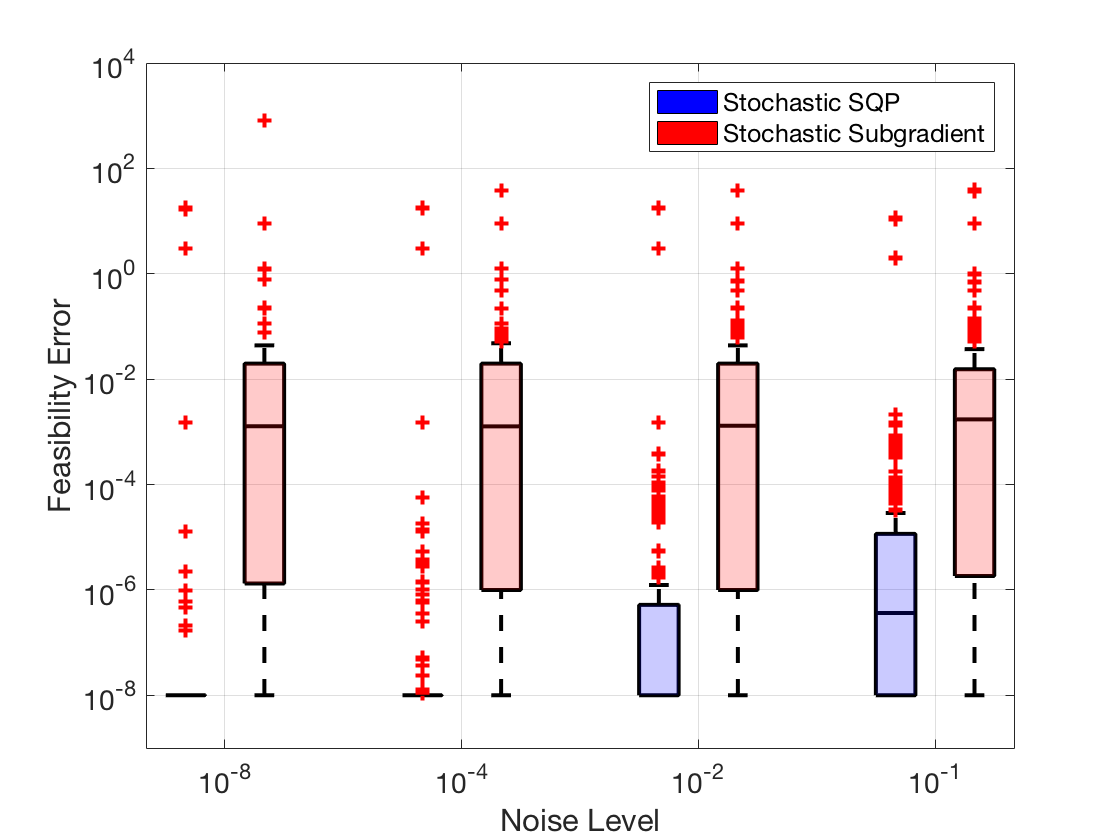}\quad 
  \includegraphics[width=0.45\textwidth,clip=true,trim=20 5 90 50]{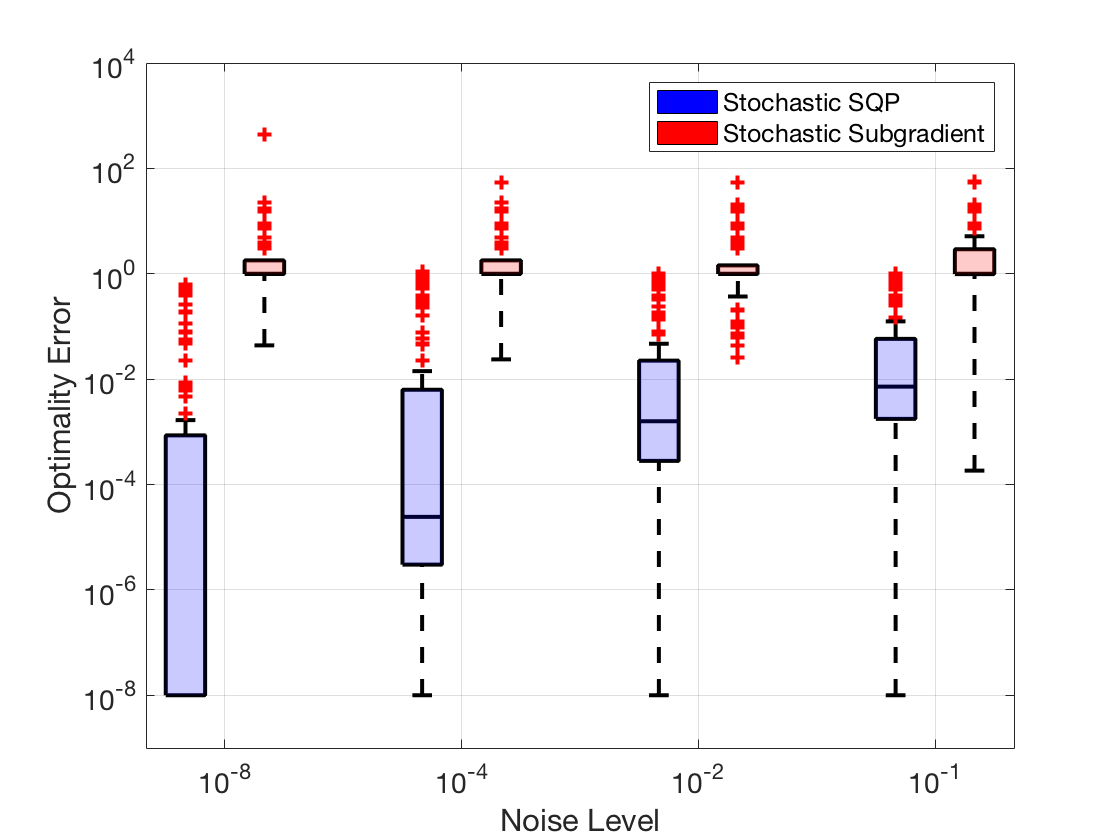}
  \caption{Box plots for feasibility errors (left) and optimality errors (right).}
  \label{fig.perf_stochastic}
\end{figure}

For each algorithm and each problem instance, we computed a resulting feasibility error and optimality error as follows.  If a run produced an iterate that was sufficiently feasible in the sense that $\|c_k\|_{\infty} \leq 10^{-6} \max\{1,\|c_0\|_{\infty}\}$ for some $k$, then, with the largest $k$ corresponding to such a feasible iterate, the feasibility error was reported as $\|c_k\|_\infty$ and the optimality error was reported as $\|g_k + J_k^Ty_k\|_\infty$, where $y_k$ was computed as a least-squares multiplier using the true gradient $g_k$ and $J_k$.  (In this manner, the optimality error is \emph{not} based on a stochastic gradient; rather, it is a true measure of optimality corresponding to the iterate $x_k$.)  On the other hand, if a run produced no sufficiently feasible iterate, then the feasibility error and optimality error were computed in this manner at the \emph{least infeasible} iterate during the run.  The results are reported in the form of box plots in Figure \ref{fig.perf_stochastic}.


Finally, let us comment on the occurrence of the event \eqref{eq.tau_small}.  In all runs of ``Stochastic SQP,'' we found that $\bar{\tau}_k \leq \tau_k^{trial}$ held $100\%$ of the time in the last 100 iterations.  In fact, for the noise levels $10^{-8}$, $10^{-4}$, $10^{-2}$, and $10^{-1}$, this inequality held in $99.92\%$, $99.10\%$, $99.22\%$, and $99.65\%$, respectively, of \emph{all} iterations.  This provides evidence that the theory offered under the event \eqref{eq.tau_small} is relevant in practice.

\section{Conclusion}\label{sec.conclusion}

We have presented, analyzed, and tested sequential quadratic optimization algorithms for solving smooth nonlinear optimization problems with equality constraints.  Our first algorithm is based on a state-of-the-art line-search SQP method, but employs a stepsize scheme based on (adaptively estimated) Lipschitz constants in place of the line search.  We have shown that this method has convergence guarantees that match those of the state-of-the-art line-search SQP method, and our numerical experiments show that the algorithm is as reliable as this state-of-the-art approach.  Based on this proposed algorithm, our second algorithm is designed to solve problems involving deterministic constraint functions, but a stochastic objective function.  We have proved that under good behavior of the merit function parameter, the algorithm possesses convergence guarantees that match those of our deterministic algorithm in expectation.  We have also argued that certain poor behavior of the merit function parameter will only occur in extreme circumstances, and other poor behavior only occurs with probability zero (and in any case can be safeguarded against).  Our numerical experiments show that our algorithm for the stochastic setting consistently and significantly outperforms a (sub)gradient method employed to minimize a penalty function, which is an algorithm that represents the current state-of-the-art in the context of \emph{stochastic} constrained optimization.

One assumption required for our analysis is that the iterates remain in an open convex set over which the objective and constraint functions and their derivatives remain bounded.  This is not ideal in the context of a stochastic algorithm, although it is more forgivable in a constrained setting than in an unconstrained setting since the algorithm is designed to be driven to the deterministic feasible region.  That being said, one could loosen this assumption if one were to apply our algorithm with $\theta = 0$.  Indeed, notice that in our analysis in \S\ref{sec.constant_tau_small}, boundedness of $\{\|g_k\|_2\}$ is primarily required in Lemma~\ref{lem.alphagd}, but with $\theta=0$ one finds directly that, for $k \geq \kbar_{\tau,\xi}$,
\bequationNN
  \E_k[\bar\alpha_k\bar\tau_kg_k^T(\dbar_k - d_k)] = \(\tfrac{\beta_k \bar\xi_{\min} \bar\tau_{\min}}{\bar\tau_{\min}L_k + \Gamma_k}\) \bar\tau_{\min} \E_k[g_k^T(\dbar_k - d_k)] = 0.
\eequationNN
Hence, our assumption about the boundedness of $\{\|g_k\|_2\}$ is only needed when $\theta > 0$.  We have proposed our algorithm for this setting since it is the context of $\theta > 0$ that allows the stepsize scheme in our algorithm to be adaptive, which has a significant benefit in terms of practical performance of the method.

\bibliographystyle{plain}
\bibliography{references}

\begin{thebibliography}{10}

\bibitem{BongConnGoulToin95}
Ingrid Bongartz, Andrew~R Conn, Nick Gould, and Ph~L Toint.
\newblock Cute: Constrained and unconstrained testing environment.
\newblock {\em ACM Transactions on Mathematical Software (TOMS)},
  21(1):123--160, 1995.

\bibitem{BottCurtNoce18}
L{\'e}on Bottou, Frank~E Curtis, and Jorge Nocedal.
\newblock Optimization methods for large-scale machine learning.
\newblock {\em SIAM Review}, 60(2):223--311, 2018.

\bibitem{ByrdHribNoce99}
R.~H. Byrd, M.~E. Hribar, and J.~Nocedal.
\newblock An interior point algorithm for large-scale nonlinear programming.
\newblock {\em SIAM Journal on Optimization}, 9(4):877--900, 1999.

\bibitem{ByrdCurtNoce08}
Richard~H. Byrd, Frank~E. Curtis, and Jorge Nocedal.
\newblock {An Inexact SQP Method for Equality Constrained Optimization}.
\newblock {\em {SIAM Journal on Optimization}}, 19(1):351--369, 2008.

\bibitem{ByrdCurtNoce10a}
Richard~H. Byrd, Frank~E. Curtis, and Jorge Nocedal.
\newblock {An Inexact Newton Method for Nonconvex Equality Constrained
  Optimization}.
\newblock {\em {Mathematical Programming}}, 122(2):273--299, 2010.

\bibitem{ByrdGilbNoce00}
Richard~H. Byrd, Jean~Charles Gilbert, and Jorge Nocedal.
\newblock A trust region method based on interior point techniques for
  nonlinear programming.
\newblock {\em Mathematical Programming}, 89(1):149--185, 2000.

\bibitem{ChenTungVeduMori18}
Changan Chen, Frederick Tung, Naveen Vedula, and Greg Mori.
\newblock {Constraint-aware deep neural network compression}.
\newblock In {\em Proceedings of the European Conference on Computer Vision
  (ECCV)}, pages 400--415, 2018.

\bibitem{ConnGoulToin92}
A.~R. Conn, N.~I.~M. Gould, and {Ph}.~L. Toint.
\newblock {\em {LANCELOT}: A Fortran package for large-scale nonlinear
  optimization (Release A)}.
\newblock Number~17 in Springer Series in Computational Mathematics. Springer
  Verlag, Heidelberg, Berlin, New York, 1992.

\bibitem{Cour43}
R.~Courant.
\newblock Variational methods for the solution of problems of equilibrium and
  vibrations.
\newblock {\em Bull. Amer. Math. Soc.}, 49(1):1--23, 01 1943.

\bibitem{CurtRobi19}
Frank~E. Curtis and Daniel~P. Robinson.
\newblock {Exploiting Negative Curvature in Deterministic and Stochastic
  Optimization}.
\newblock {\em {Mathematical Programming, Series B}}, 176(1):69--94, 2019.

\bibitem{DolaMore02}
Elizabeth~D Dolan and Jorge~J Mor{\'e}.
\newblock Benchmarking optimization software with performance profiles.
\newblock {\em Mathematical programming}, 91(2):201--213, 2002.

\bibitem{Flet87}
R.~Fletcher.
\newblock {\em {Practical Methods of Optimization}}.
\newblock John Wiley and Sons, Chichester, UK, 2 edition, 1987.

\bibitem{goldfarb2017linear}
Donald Goldfarb, Garud Iyengar, and Chaoxu Zhou.
\newblock Linear convergence of stochastic {F}rank-{W}olfe variants.
\newblock {\em arXiv preprint arXiv:1703.07269}, 2017.

\bibitem{Han77}
S.~P. Han.
\newblock {A Globally Convergent Method for Nonlinear Programming}.
\newblock {\em Journal of Optimization Theory and Applications},
  22(3):297--309, 1977.

\bibitem{HanMang79}
S.~P. Han and O.~L. Mangasarian.
\newblock {Exact Penalty Functions in Nonlinear Programming}.
\newblock {\em Mathematical Programming}, 17:251--269, 1979.

\bibitem{hazan2016variance}
Elad Hazan and Haipeng Luo.
\newblock Variance-reduced and projection-free stochastic optimization.
\newblock In {\em International Conference on Machine Learning}, pages
  1263--1271, 2016.

\bibitem{Hest69}
M.~R. Hestenes.
\newblock {Multiplier and Gradient Methods}.
\newblock {\em Journal of Optimization Theory and Applications}, 4:303--320,
  1969.

\bibitem{KumaSoumMhamHara18}
Soumava Kumar~Roy, Zakaria Mhammedi, and Mehrtash Harandi.
\newblock {Geometry aware constrained optimization techniques for deep
  learning}.
\newblock In {\em Proceedings of the IEEE Conference on Computer Vision and
  Pattern Recognition}, pages 4460--4469, 2018.

\bibitem{locatello2019stochastic}
Francesco Locatello, Alp Yurtsever, Olivier Fercoq, and Volkan Cevher.
\newblock Stochastic {F}rank-{W}olfe for composite convex minimization.
\newblock In {\em Advances in Neural Information Processing Systems}, pages
  14269--14279, 2019.

\bibitem{lu2020generalized}
Haihao Lu and Robert~M Freund.
\newblock Generalized stochastic {F}rank-{W}olfe algorithm with stochastic
  “substitute” gradient for structured convex optimization.
\newblock {\em Mathematical Programming}, pages 1--33, 2020.

\bibitem{NandPathAbhiSing19}
Yatin Nandwani, Abhishek Pathak, and Parag Singla.
\newblock {A primal-dual formulation for deep learning with constraints}.
\newblock In {\em Advances in Neural Information Processing Systems}, pages
  12157--12168, 2019.

\bibitem{negiar2020stochastic}
Geoffrey N{\'e}giar, Gideon Dresdner, Alicia Tsai, Laurent~El Ghaoui, Francesco
  Locatello, and Fabian Pedregosa.
\newblock Stochastic {F}rank-{W}olfe for constrained finite-sum minimization.
\newblock {\em arXiv preprint arXiv:2002.11860}, 2020.

\bibitem{Nest04}
Yurii Nesterov.
\newblock {\em Introductory Lectures on Convex Optimization}.
\newblock Applied Optimization. Springer Science+Business Media New York, 2004.

\bibitem{NoceWrig06}
Jorge Nocedal and Stephen Wright.
\newblock {\em Numerical optimization}.
\newblock Springer Series in Operations Research and Financial Engineering.
  Springer-Verlag New York, 2006.

\bibitem{Powe69}
M.~J.~D. Powell.
\newblock {A Method for Nonlinear Constraints in Minimization Problems}.
\newblock In R.~Fletcher, editor, {\em Optimization}, pages 283--298. Academic
  Press, London and New York, 1969.

\bibitem{Powe78}
M.~J.~D. Powell.
\newblock {A Fast Algorithm for Nonlinearly Constrained Optimization
  Calculations}.
\newblock In {\em Numerical Analysis}, Lecture Notes in Mathematics, pages
  144--157. Springer, Berlin, Heidelberg, 1978.

\bibitem{RaviDinhLokhSing19}
Sathya~N Ravi, Tuan Dinh, Vishnu~Suresh Lokhande, and Vikas Singh.
\newblock {Explicitly imposing constraints in deep networks via conditional
  gradients gives improved generalization and faster convergence}.
\newblock In {\em Proceedings of the AAAI Conference on Artificial
  Intelligence}, volume~33, pages 4772--4779, 2019.

\bibitem{reddi2016stochastic}
Sashank~J Reddi, Suvrit Sra, Barnab{\'a}s P{\'o}czos, and Alex Smola.
\newblock Stochastic {F}rank-{W}olfe methods for nonconvex optimization.
\newblock In {\em 2016 54th Annual Allerton Conference on Communication,
  Control, and Computing (Allerton)}, pages 1244--1251. IEEE, 2016.

\bibitem{ShapDentRusz09}
A.~Shapiro, D.~Dentcheva, and A.~Ruszczy{\'n}ski.
\newblock {\em Lectures on Stochastic Programming: Modeling and Theory}.
\newblock SIAM, 2009.

\bibitem{Tong12}
Yung~Liang Tong.
\newblock {\em The multivariate normal distribution}.
\newblock Springer Science \& Business Media, 2012.

\bibitem{WaecBieg06}
A.~Waechter and L.~T. Biegler.
\newblock {On the implementation of an interior-point filter line-search
  algorithm for large-scale nonlinear programming}.
\newblock {\em Mathematical Programming}, 106:25--57, 2006.

\bibitem{Wils63}
R.~B. Wilson.
\newblock {\em {A Simplicial Algorithm for Concave Programming}}.
\newblock Ph.D. Thesis, Graduate School of Business Administration, Harvard
  University, Cambridge, MA, USA, 1963.

\bibitem{zhang2020one}
Mingrui Zhang, Zebang Shen, Aryan Mokhtari, Hamed Hassani, and Amin Karbasi.
\newblock One sample stochastic {F}rank-{W}olfe.
\newblock In {\em International Conference on Artificial Intelligence and
  Statistics}, pages 4012--4023, 2020.

\end{thebibliography}

\end{document}